\documentclass[a4paper,fleqn]{cas-sc}

\usepackage[numbers]{natbib}
\usepackage{booktabs}
\usepackage{siunitx}

\def\tsc#1{\csdef{#1}{\textsc{\lowercase{#1}}\xspace}}
\tsc{WGM}
\tsc{QE}

\begin{document}
\let\WriteBookmarks\relax
\def\floatpagepagefraction{1}
\def\textpagefraction{.001}

\shorttitle{}    

\shortauthors{}  

\title [mode = title]{
Domain-decomposed parallelization of 
B-spline based s-version of the finite element method 
via generalized graph abstraction
}

\author[1]{Nozomi Magome}
\ead{magome.nozomi.sw@alumni.tsukuba.ac.jp}
\credit{Data curation, Formal analysis, Investigation, Methodology,
Software, Validation, Visualization, Funding acquisition,
Writing -- original draft}

\author[2]{Naoki Morita}
\ead{nmorita@kz.tsukuba.ac.jp}
\credit{Methodology, Software, Funding acquisition,
Writing -- review and editing, Supervision}

\author[3]{Shigeki Kaneko}
\ead{kaneko.shigeki@nitech.ac.jp}
\credit{Methodology, Writing -- review and editing, Supervision}

\author[2]{Naoto Mitsume}
\cormark[1]
\ead{mitsume@kz.tsukuba.ac.jp}
\credit{Conceptualization, Project administration, Funding acquisition,
Writing -- review and editing, Supervision, Resources}

\affiliation[1]{
  organization={Degree Programs in Systems and Information Engineering,
                University of Tsukuba},
  addressline={Tennodai 1-1-1},
  city={Tsukuba},
  postcode={305-8573},
  state={Ibaraki},
  country={Japan}
}

\affiliation[2]{
  organization={Institute of Systems and Information Engineering,
                University of Tsukuba},
  addressline={Tennodai 1-1-1},
  city={Tsukuba},
  postcode={305-8573},
  state={Ibaraki},
  country={Japan}
}

\affiliation[3]{
  organization={Graduate School of Engineering,
                Nagoya Institute of Technology},
  addressline={Gokiso-cho, Showa-ku},
  city={Nagoya},
  postcode={466-8555},
  state={Aichi},
  country={Japan}
}

\cortext[1]{Corresponding author}

\begin{abstract}
The s-version of the finite element method (SFEM) enables locally high-resolution analysis by superimposing independently defined finite element meshes. However, domain-decomposed parallelization is nontrivial because complex interactions arise among degrees of freedom distributed over multiple meshes. In this study, we propose a method for constructing a graph structure that uniformly represents interactions among computational points, including both intra- and inter-mesh interactions, based on the overlap of basis-function supports. We apply the proposed graph representation to the B-spline based SFEM (BSFEM), a high-accuracy SFEM formulation previously proposed by the authors. The resulting graph partition enables the consistent assignment of degrees of freedom and elements to processes and the construction of the MPI communication structure, thereby realizing domain-decomposition-based distributed-memory parallelization of BSFEM. To the best of the authors’ knowledge, this BSFEM implementation constitutes the first domain-decomposition-based distributed-memory parallelization of an SFEM-based method. Furthermore, as an example demonstrating the utility of the proposed graph representation, we apply cost-weighted graph partitioning in which the matrix-generation costs specific to BSFEM are incorporated into node weights, and demonstrate effective static load balancing that accounts for the nonuniform matrix-generation workload.
\end{abstract}


\begin{highlights}
\item A novel generalized graph representation for BSFEM is proposed.
\item The graph captures intra- and inter-mesh interactions.
\item The graph enables the first domain-decomposed parallelization of an SFEM.
\item Cost weighting shows the graph's extensibility and improves static load balancing.
\end{highlights}

\begin{keywords}
s-version of the finite element method\sep
B-spline basis functions\sep
Distributed-memory parallelization\sep
Domain decomposition\sep
Overlapping meshes\sep
Load balancing
\end{keywords}

\maketitle

\section{Introduction}
In engineering simulations, it is often necessary to resolve localized regions with steep solution gradients, such as crack tips, localized heat sources, and boundary layers, at high spatial resolution while limiting computational cost. In conventional finite difference methods and finite element methods (FEM), a common approach is to locally refine a single mesh conforming to the problem geometry \citep{thames1977numerical,hughes2005isogeometric,el2022anisotropic}. However, for problems involving complex geometries or regions requiring high resolution that move or evolve, mesh refinement and regeneration can be computationally expensive, and element distortion may degrade numerical accuracy and stability \citep{mittal2005immersed}.
One approach to this problem is to superimpose independently generated meshes, as in overset grid methods \citep{benek1983flexible} and composite grid methods \citep{steger1987use,chesshire1990composite}, which have been applied to fluid, solid, and multiphysics problems \citep{tang2003overset,appelo2012numerical,henshaw2009composite}. Such multiple-mesh approaches facilitate mesh generation for complex geometries while helping maintain mesh quality.

Within the FEM framework, the s-version of the finite element method (SFEM) has been proposed as a related multiple-mesh approach \citep{fish1992s,fish1993multiscale,FISH1994135}. In SFEM, a global mesh representing the entire computational domain and a local mesh covering a region of interest are defined independently, and the solution is represented by superimposing the basis functions of the two meshes. This allows high spatial resolution to be introduced only where needed without regenerating the global mesh. Whereas overset methods generally require inter-mesh interpolation or hole cutting, SFEM incorporates the basis functions of both meshes into the same variational formulation and couples them directly within a single finite element discretization.
SFEM has been applied to various problems requiring locally high resolution, including fracture and damage analyses \citep{lee2004combined,okada2005fracture,nakasumi2008crack,kikuchi2014fatigue,xu2018study,kishi2020dynamic}, multiscale stress analyses \citep{sun2018variant,sakata2022mesh}, and shape optimization and elasticity analyses \citep{wang2006moving,yue2005adaptive}.
Conventional SFEM has two major issues: reduced numerical integration accuracy in the overlapping region and singularity of the coefficient matrix caused by linear dependence between basis functions defined on different meshes \citep{FISH1993363,he2023strategy,fan2008rs,ooya2009linear}. To address these issues, the authors proposed the B-spline based s-version of the finite element method (BSFEM), which employs Lagrange basis functions on the local mesh and highly continuous B-spline basis functions on the global mesh \citep{magome2024higher}. Numerical analyses of steady Poisson problems demonstrated that BSFEM improves numerical integration accuracy while avoiding coefficient-matrix singularity and loss of positive definiteness.

Distributed-memory parallel computing is required to apply SFEM to large-scale three-dimensional problems. In conventional FEM, the mesh connectivity can be represented as a graph, and the computational domain can be decomposed using graph partitioning libraries such as METIS to distribute the computational workload while limiting interprocess communication \citep{karypis1998fast}. In SFEM, however, interactions between degrees of freedom arise not only within each mesh but also between different meshes in the overlapping region. Therefore, when the meshes are partitioned independently, the data dependencies induced by inter-mesh coupling are not directly reflected in the domain decomposition.
Embarrassingly parallel computations based on SFEM, in which mutually independent analyses are executed in parallel, have been reported \citep{suwa2023parallel}. However, to the authors' knowledge, no method has been reported for decomposing a single SFEM discretization and performing distributed-memory parallel computation while accounting for both intra- and inter-mesh interactions between degrees of freedom. Its realization therefore requires a unified representation of the interactions between degrees of freedom across multiple independent meshes that can be used for domain decomposition.

In this study, we propose a parallelization framework that represents the interactions between basis functions in BSFEM as a single graph. Each computational point associated with a basis function is represented as a graph node, and an edge is defined between graph nodes for which a coefficient-matrix entry can arise from the discretization. This representation expresses intra- and inter-mesh interactions within a single graph without distinguishing between them. Because the definition is independent of the type of basis function, interactions between the B-spline and Lagrange bases used in BSFEM can also be treated within the same framework. This generalization enables existing graph partitioning libraries \citep{karypis1998fast} to be applied to domain decomposition in BSFEM. The resulting graph partition is used to construct the ownership of degrees of freedom and elements by each process, together with the MPI communication structure, thereby enabling domain-decomposition-based distributed-memory parallel computation of BSFEM. To the authors' knowledge, this is the first domain-decomposition-based distributed-memory parallel computation for an SFEM-based method.
Furthermore, the parallel implementation revealed that the computational cost of matrix generation in BSFEM is spatially nonuniform. We therefore introduce static load balancing based on weighted graph partitioning, in which the matrix-generation cost is represented by graph-node weights. The proposed graph representation thus allows both the interactions between computational points and the nonuniform computational workload to be incorporated into the domain decomposition, enabling practically useful parallel performance.

The remainder of this paper is organized as follows. Section~\ref{sec:poisson_formulation_SFEM} presents the formulations of SFEM and BSFEM. Section~\ref{sec:parallelization} describes the proposed graph-based parallelization framework. Section~\ref{sec:verification_parallelization} presents the verification results and parallel performance of the proposed method. Finally, Section~\ref{sec:conclusions} gives the conclusions.

\section{Formulation of s-version of the finite element method (SFEM) and B-spline based SFEM (BSFEM)}
\label{sec:poisson_formulation_SFEM}
\subsection{SFEM discretization of the Poisson equation}
\label{subsec:sfem_discretization}

This subsection defines the discretized matrix system of SFEM, which is the target of the parallelization method proposed in this study.
The general formulation of SFEM and the details of the Galerkin discretization follow our previous work \cite{magome2024higher}.
Therefore, this subsection presents only the governing equation, the representation of the solution in SFEM, and the resulting block matrix structure required in this paper.

As a fundamental verification problem, this study considers the steady Poisson equation with Dirichlet boundary conditions imposed on the entire boundary:
\begin{subequations}
\label{eq:poisson}
\begin{align}
    \Delta u + f &= 0 
    \quad \mathrm{in} \; \Omega, \\
    u &= g 
    \quad \mathrm{on} \; \Gamma_{\mathrm{D}},
\end{align}
\end{subequations}
where $\Omega$ is the analysis domain, $\Gamma_{\mathrm{D}}=\partial\Omega$ is the Dirichlet boundary, $u$ is the trial solution, and $f$ and $g$ are prescribed functions.
The standard weak form is given as the problem of finding $u$ such that
\begin{equation}
    a_{\Omega}(w,u)=L_{\Omega}(w)
\end{equation}
for all test functions $w$.
Here,
\begin{align}
    a_{\omega}(w,u) &= \int_{\omega} \nabla w \cdot \nabla u \, d\Omega, \\
    L_{\omega}(w) &= \int_{\omega} w f \, d\Omega.
\end{align}

In SFEM \cite{fish1992s}, a local mesh is superimposed on a subdomain where higher resolution is required, while a global mesh covers the entire analysis domain.
In this study, the global domain is defined as $\Omega^{\mathrm{G}}=\Omega$, and the local domain $\Omega^{\mathrm{L}}$ is assumed to satisfy
\begin{equation}
    \Omega^{\mathrm{L}} \subseteq \Omega^{\mathrm{G}}.
\end{equation}
However, the boundary of $\Omega^{\mathrm{L}}$ is allowed to partially coincide with $\partial\Omega^{\mathrm{G}}$.
The global and local meshes are not required to be conforming.

\begin{figure}
    \centering
    \includegraphics[bb=0 0 993 564,width=9cm]{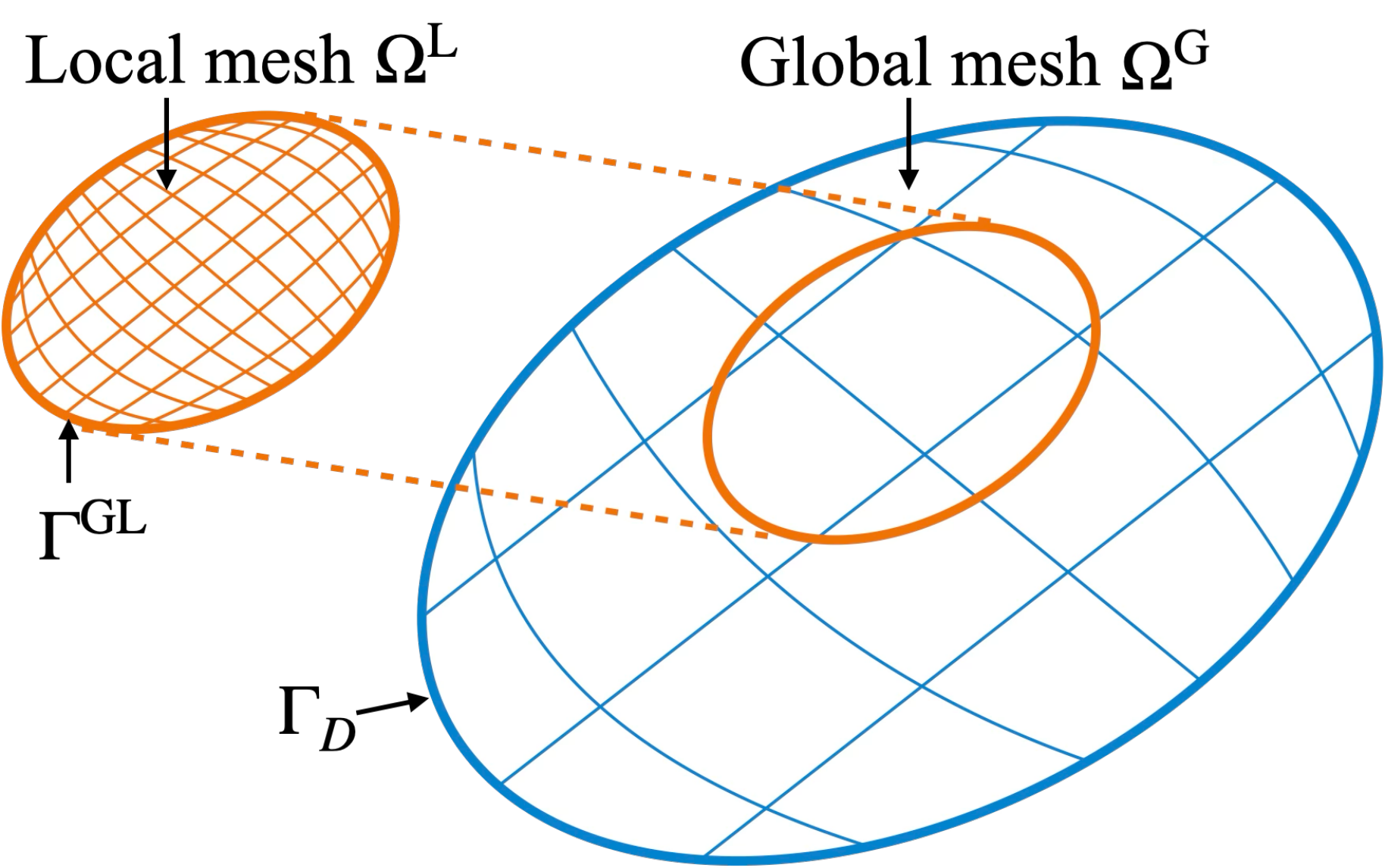}
    \caption{Global and local meshes in the s-version of the finite element method. The global mesh covers the entire domain $\Omega^{\mathrm{G}}$, while the local mesh is independently overlaid on the subdomain $\Omega^{\mathrm{L}}$.}
    \label{fig:sfem_meshes}
\end{figure}

In SFEM, the trial solution is represented as the sum of a global component $u^{\mathrm{G}}$ and a local component $u^{\mathrm{L}}$.
That is,
\begin{equation}
    u =
    \begin{cases}
        u^{\mathrm{G}} 
        & \mathrm{in} \; \Omega^{\mathrm{G}} \setminus \Omega^{\mathrm{L}},\\
        u^{\mathrm{G}} + u^{\mathrm{L}}
        & \mathrm{in} \; \Omega^{\mathrm{L}}.
    \end{cases}
    \label{eq:sfem_solution}
\end{equation}
The Dirichlet boundary condition is imposed on the global solution along the external boundary as
\begin{equation}
    u^{\mathrm{G}} = g
    \quad \mathrm{on} \; \Gamma_{\mathrm{D}},
\end{equation}
whereas the following condition is imposed on the local solution along the boundary of the local domain, $\Gamma^{\mathrm{GL}}=\partial\Omega^{\mathrm{L}}$:
\begin{equation}
    u^{\mathrm{L}} = 0
    \quad \mathrm{on} \; \Gamma^{\mathrm{GL}}.
    \label{eq:local_dirichlet}
\end{equation}
This condition eliminates the local component on $\Gamma^{\mathrm{GL}}$ and preserves continuity with the global solution.

Corresponding to Eq.~\eqref{eq:sfem_solution}, the test function is also represented as
\begin{equation}
    w =
    \begin{cases}
        w^{\mathrm{G}} 
        & \mathrm{in} \; \Omega^{\mathrm{G}} \setminus \Omega^{\mathrm{L}},\\
        w^{\mathrm{G}} + w^{\mathrm{L}}
        & \mathrm{in} \; \Omega^{\mathrm{L}}.
    \end{cases}
\end{equation}
The weak form in SFEM is then written as
\begin{equation}
    a_{\Omega}'(w,u)=L_{\Omega}'(w),
    \label{eq:sfem_weak_form}
\end{equation}
where
\begin{align}
    a_{\Omega}'(w,u)
    &=
    a_{\Omega^{\mathrm{G}}}(w^{\mathrm{G}},u^{\mathrm{G}})
    + a_{\Omega^{\mathrm{L}}}(w^{\mathrm{G}},u^{\mathrm{L}})
    + a_{\Omega^{\mathrm{L}}}(w^{\mathrm{L}},u^{\mathrm{G}})
    + a_{\Omega^{\mathrm{L}}}(w^{\mathrm{L}},u^{\mathrm{L}}), \\
    L_{\Omega}'(w)
    &=
    L_{\Omega^{\mathrm{G}}}(w^{\mathrm{G}})
    + L_{\Omega^{\mathrm{L}}}(w^{\mathrm{L}}).
\end{align}

Next, let the global and local basis functions be denoted by
$N_j^{\mathrm{G}}$ and $N_k^{\mathrm{L}}$, respectively.
Let $\mathcal{I}^{\mathrm{G}}_0$ be the index set of the global nodes excluding those on the Dirichlet boundary, and let $\mathcal{I}^{\mathrm{L}}$ be the index set of the local nodes.
If $g^h$ denotes a function satisfying the Dirichlet boundary condition, the trial solution is approximated as
\begin{align}
    (u^{\mathrm{G}})^h
    &=
    \sum_{j \in \mathcal{I}^{\mathrm{G}}_0}
    N_j^{\mathrm{G}} q_j^{\mathrm{G}}
    + g^h, \\
    (u^{\mathrm{L}})^h
    &=
    \sum_{k \in \mathcal{I}^{\mathrm{L}}}
    N_k^{\mathrm{L}} q_k^{\mathrm{L}}.
\end{align}

This Galerkin discretization yields the following system of linear equations:
\begin{equation}
    \boldsymbol{K}\boldsymbol{q}
    =
    \boldsymbol{F},
    \label{eq:system_sfem}
\end{equation}
where
\begin{equation}
    \boldsymbol{K}
    =
    \begin{bmatrix}
        \boldsymbol{K}^{\mathrm{GG}} & \boldsymbol{K}^{\mathrm{GL}} \\
        \boldsymbol{K}^{\mathrm{LG}} & \boldsymbol{K}^{\mathrm{LL}}
    \end{bmatrix},
    \quad
    \boldsymbol{q}
    =
    \begin{bmatrix}
        \boldsymbol{q}^{\mathrm{G}} \\
        \boldsymbol{q}^{\mathrm{L}}
    \end{bmatrix},
    \quad
    \boldsymbol{F}
    =
    \begin{bmatrix}
        \boldsymbol{F}^{\mathrm{G}} \\
        \boldsymbol{F}^{\mathrm{L}}
    \end{bmatrix}.
    \label{eq:sfem_block_matrix}
\end{equation}

The components of each submatrix are defined, for $i,j \in \mathcal{I}^{\mathrm{G}}_0$ and
$k,l \in \mathcal{I}^{\mathrm{L}}$, as
\begin{align}
    K_{ij}^{\mathrm{GG}}
    &=
    a_{\Omega^{\mathrm{G}}}
    \left(
        N_i^{\mathrm{G}},
        N_j^{\mathrm{G}}
    \right), \\
    K_{ik}^{\mathrm{GL}}
    &=
    a_{\Omega^{\mathrm{L}}}
    \left(
        N_i^{\mathrm{G}},
        N_k^{\mathrm{L}}
    \right), \\
    K_{lj}^{\mathrm{LG}}
    &=
    a_{\Omega^{\mathrm{L}}}
    \left(
        N_l^{\mathrm{L}},
        N_j^{\mathrm{G}}
    \right), \\
    K_{lk}^{\mathrm{LL}}
    &=
    a_{\Omega^{\mathrm{L}}}
    \left(
        N_l^{\mathrm{L}},
        N_k^{\mathrm{L}}
    \right).
\end{align}
The components of the right-hand-side vector are given by
\begin{align}
    F_i^{\mathrm{G}}
    &=
    L_{\Omega^{\mathrm{G}}}
    \left(
        N_i^{\mathrm{G}}
    \right)
    -
    a_{\Omega^{\mathrm{G}}}
    \left(
        N_i^{\mathrm{G}},
        g^h
    \right), \\
    F_l^{\mathrm{L}}
    &=
    L_{\Omega^{\mathrm{L}}}
    \left(
        N_l^{\mathrm{L}}
    \right)
    -
    a_{\Omega^{\mathrm{L}}}
    \left(
        N_l^{\mathrm{L}},
        g^h
    \right).
\end{align}
The detailed derivation of the above equations is given in our previous work \cite{magome2024higher}.
Here, $\boldsymbol{K}^{\mathrm{GG}}$ represents the interactions within the global mesh, whereas $\boldsymbol{K}^{\mathrm{LL}}$ represents the interactions within the local mesh.
On the other hand, $\boldsymbol{K}^{\mathrm{GL}}$ and $\boldsymbol{K}^{\mathrm{LG}}$ represent the interactions in the overlapping region between the global and local meshes.
Hereafter, these matrices are referred to as coupling matrices.
The parallelization method developed in this study targets the SFEM-specific block matrix structure shown in Eq.~\eqref{eq:sfem_block_matrix}.

\subsection{Construction of basis functions}
\label{sq:BSFEM}

This subsection describes the construction of the Lagrange and B-spline basis functions used in this study.
The general definitions of both basis functions and their basic use in the B-spline based s-version of the finite element method (BSFEM) follow our previous work \cite{magome2024higher}.
Therefore, this subsection focuses on the analysis conditions required for reproducibility and on the element-local construction of B-spline basis functions using Bézier extraction introduced in this study.
For simplicity, the following description is given for the one-dimensional case.
Multidimensional parametric basis functions are defined as tensor products of the one-dimensional basis functions.

The Lagrange basis functions are defined in terms of the reference coordinate $\hat{\xi}$ on each element, as in the standard finite element method.
The $p$th-order Lagrange basis function is defined on the parent element $\hat{\Omega}_e=[-1,1]$ with nodes
$\hat{\xi}_i \, (i=1,\ldots,p+1)$ as
\begin{equation}
    l_i^p(\hat{\xi})
    =
    \prod_{\substack{j=1 \\ j \neq i}}^{p+1}
    \frac{\hat{\xi}-\hat{\xi}_j}{\hat{\xi}_i-\hat{\xi}_j}.
    \label{eq:lagrange}
\end{equation}
This basis function satisfies the Kronecker delta property at the nodes and is $C^0$-continuous across element boundaries.

The B-spline basis functions are defined by a nondecreasing knot vector
\begin{equation}
    \boldsymbol{\Xi}
    =
    \left\{
    \xi_1,\xi_2,\ldots,\xi_{n_{\mathrm{knot}}}
    \right\}^{\mathrm{T}}.
    \label{eq:knot_vec}
\end{equation}
Here, $n_{\mathrm{knot}}=n_{\mathrm{cp}}+p+1$, where
$n_{\mathrm{cp}}$ is the number of control points and $p$ is the degree of the B-spline basis functions.

In this study, an open knot vector is used so that the basis functions are interpolatory at both ends.
The multiplicity of all internal knots is set to one, resulting in $C^{p-1}$ continuity of the B-spline basis functions across internal element boundaries.
The geometry is represented by a linear combination of the B-spline basis functions and control points.
However, because the control variables do not generally coincide with the physical field values, their coordinates are determined using the mesh generation method of Otoguro et al. \cite{otoguro2017space}.
The Dirichlet boundary condition is imposed by exploiting the interpolatory property at the endpoints of the open knot vector and fixing the values of the control points corresponding to the Dirichlet boundary to the prescribed values.

In this study, Bézier extraction \cite{Borden2011} is used for the element-wise construction of the B-spline basis functions.
Bézier extraction represents B-spline basis functions as linear combinations of Bernstein polynomials defined on each element.
This enables basis function evaluation on each element to be performed as a local matrix operation using only the Bézier extraction operator associated with that element.

Let $\mathbf{M}_p^e(\hat{\xi})$ denote the vector of B-spline element basis functions of degree $p$ on element $e$, and let
$\mathbf{B}_p(\hat{\xi})$ denote the vector of Bernstein polynomials on the parent element $\hat{\Omega}_e=[-1,1]$.
Then,
\begin{equation}
    \mathbf{M}_p^e(\hat{\xi})
    =
    \mathbf{C}_p^e
    \mathbf{B}_p(\hat{\xi})
    \label{eq:bezier_extraction}
\end{equation}
is obtained.
Here, $\mathbf{C}_p^e$ is the Bézier extraction operator corresponding to element $e$.
The vectors are given by
\begin{align}
    \mathbf{M}_p^e(\hat{\xi})
    &=
    \left\{
    M_{1,p}^e(\hat{\xi}),
    M_{2,p}^e(\hat{\xi}),
    \ldots,
    M_{p+1,p}^e(\hat{\xi})
    \right\}^{\mathrm{T}},
    \\
    \mathbf{B}_p(\hat{\xi})
    &=
    \left\{
    B_{1,p}(\hat{\xi}),
    B_{2,p}(\hat{\xi}),
    \ldots,
    B_{p+1,p}(\hat{\xi})
    \right\}^{\mathrm{T}}.
\end{align}
The Bernstein polynomials are defined as
\begin{equation}
    B_{i,p}(\hat{\xi})
    =
    \binom{p}{i-1}
    \left(
        \frac{1-\hat{\xi}}{2}
    \right)^{p-i+1}
    \left(
        \frac{1+\hat{\xi}}{2}
    \right)^{i-1},
    \quad
    i=1,2,\ldots,p+1.
    \label{eq:bernstein_basis}
\end{equation}

The computation of the stiffness matrix requires the derivatives of the basis functions.
When Bézier extraction is used, the derivatives of the B-spline basis functions can also be evaluated locally as
\begin{equation}
    \frac{d\mathbf{M}_p^e}{d\hat{\xi}}
    =
    \mathbf{C}_p^e
    \frac{d\mathbf{B}_p}{d\hat{\xi}}.
    \label{eq:bezier_extraction_derivative}
\end{equation}
Therefore, the Cox--de Boor recursion formula does not need to be evaluated directly at each integration point.
In the conventional evaluation based on the Cox--de Boor recursion formula, the basis functions must be constructed by referring to the knot information that determines their support.
In contrast, in Bézier extraction, the required knot information is condensed into the element-wise operators $\mathbf{C}_p^e$ during preprocessing.
Thus, when computing element stiffness matrices, each process can evaluate the B-spline basis functions using only the Bernstein polynomials and Bézier extraction operators associated with its assigned elements.
This property allows the evaluation of B-spline basis functions to be treated as an element-local operation similar to that in the standard finite element method, thereby facilitating a domain-decomposed parallel implementation.

\subsection{Computational algorithm for BSFEM}
\label{sq:BSFEM_algorithm}

This subsection describes the computational algorithm used in the numerical implementation of BSFEM.
The implementation in this study consists of the following four phases.
The main focus here is on the data structures and computational procedures required for the generation of the SFEM-specific matrices.

\begin{enumerate}
\item Mesh generation

As a preprocessing step, the global and local meshes are generated independently.
In this study, assuming application to local high-resolution analysis in interface-capturing methods, B-spline basis functions defined on an orthogonal structured grid are used for the global mesh.
In contrast, Lagrange basis functions are used for the local mesh.
Furthermore, Bézier extraction \cite{Borden2011} is applied to the B-spline basis functions so that they can be evaluated on each element using a data structure analogous to that of the standard finite element method.

\item Construction of the geometric database and coupling map

Because this study considers analyses with fixed meshes, geometric information that remains unchanged during the analysis is computed before matrix generation and stored as a database.
Specifically, for each element of the global and local meshes, the Jacobian determinant, basis function values, and physical derivatives at the integration points are precomputed.

For the generation of the coupling matrices, local elements are used as the integration domains.
Therefore, at each integration point $\boldsymbol{x}_{\mathrm{ip}}^{\mathrm{L}}$ of a local element $e^{\mathrm{L}}$, the corresponding global basis function values and derivatives must be evaluated.
To this end, the global element containing $\boldsymbol{x}_{\mathrm{ip}}^{\mathrm{L}}$ and the reference coordinates within that element are determined in advance.
In this study, by exploiting the fact that the global mesh is an orthogonal structured grid, the global element corresponding to each local integration point is identified using a bucket search.
This enables the global element search to be performed in constant time without using a tree-based search.

The obtained global element index, reference coordinates, and the global basis function values and derivatives at that point are stored in association with the corresponding local integration point.
In this study, this correspondence is referred to as the coupling map.

\item Generation of submatrices

Using the precomputed geometric database and coupling map,
$\boldsymbol{K}^{\mathrm{GG}}$,
$\boldsymbol{K}^{\mathrm{LL}}$,
$\boldsymbol{K}^{\mathrm{GL}}$,
and
$\boldsymbol{K}^{\mathrm{LG}}$
are generated.

The submatrix $\boldsymbol{K}^{\mathrm{GG}}$ associated with the global mesh is generated using standard finite element procedures with B-spline basis functions.
The submatrix $\boldsymbol{K}^{\mathrm{LL}}$ associated with the local mesh is generated using standard finite element procedures with Lagrange basis functions.
The coupling matrices $\boldsymbol{K}^{\mathrm{GL}}$ and
$\boldsymbol{K}^{\mathrm{LG}}$ are generated by evaluating the interactions between the local basis functions and the corresponding global basis functions at the integration points of the local elements.
At this stage, the corresponding global elements and basis function values are stored in the coupling map; therefore, no geometric search is performed during matrix generation.

\item Linear solver and post-processing

Finally, the degrees of freedom corresponding to the Dirichlet boundary conditions are fixed, and the resulting linear system is solved using a linear solver.
The obtained global and local solution vectors are combined according to the definition of the SFEM solution in Eq.~\eqref{eq:sfem_solution} and then output as physical quantities.

\end{enumerate}

\section{Graph-based domain-decomposition parallelization algorithm}
\label{sec:parallelization}
In this study, we propose a graph-based domain decomposition method for the B-spline based s-version of the finite element method (BSFEM), which involves complex mesh structures.
In this section, we first define a nodal graph based on the overlap of basis functions, and then describe its domain decomposition and the parallel linear algebra operations based on the resulting graph structure.

\subsection{Definition of the nodal graph}
\label{sq:nodal_graph}

\begin{figure*}
    \centering
    \includegraphics[bb=0 0 1921 467,width=13cm]{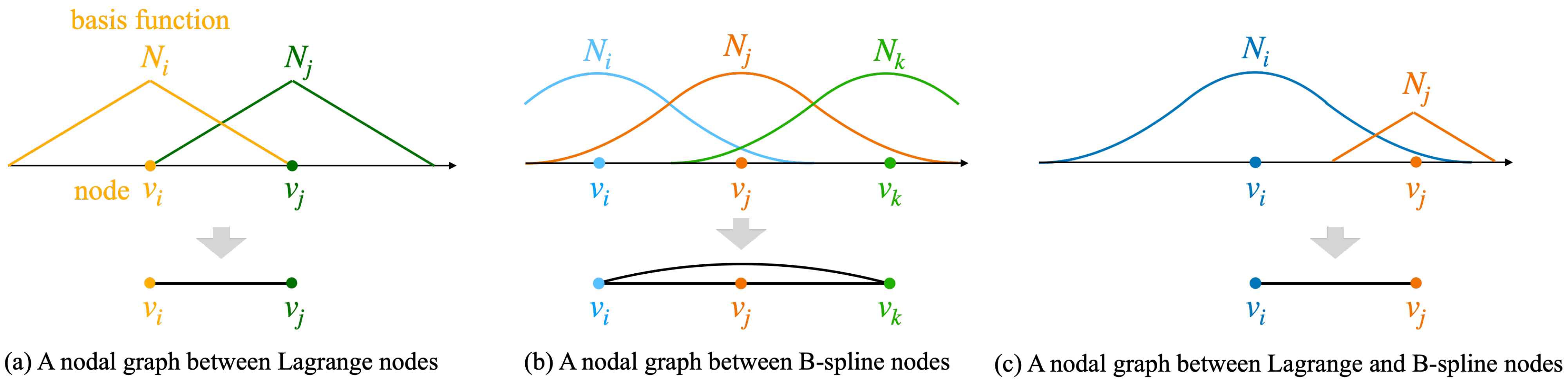}
    \caption{Construction of nodal graphs based on the overlap of basis functions.
    An edge is defined between two nodes when the corresponding basis functions are nonzero at the same point in the domain.
    This definition provides a unified graph representation for
    (a) Lagrange basis functions,
    (b) B-spline basis functions, and
    (c) superimposed meshes containing both Lagrange and B-spline nodes.}
    \label{fig:node_graph}
\end{figure*}

In this study, the connectivity among computational points on finite element meshes is represented as a graph, which is referred to as the nodal graph.
The nodal graph is defined as an undirected graph and includes self-edges.
Both nodes associated with Lagrange basis functions and control points associated with B-spline basis functions are treated as computational points corresponding to basis functions.
Hereafter, these computational points are collectively referred to as nodes.
In the nodal graph, these nodes are regarded as graph nodes, and the node set is denoted by
$V=\{v_1,v_2,\dots,v_{|V|}\}$, where $|V|$ is the total number of nodes.

Let $N_i(\boldsymbol{x})$ be the basis function corresponding to node $v_i$.
In this study, if two basis functions $N_i$ and $N_j$ are simultaneously nonzero at some point in the domain $\Omega$, the corresponding nodes $v_i$ and $v_j$ are connected by an edge.
Thus, the nodal graph is defined as
\begin{equation}
    G^{\mathrm{node}}
    =
    (V,E),
    \label{eq:nodal_graph}
\end{equation}
and its edge set $E$ is defined by
\begin{equation}
    (v_i,v_j)\in E
    \quad \Longleftrightarrow \quad
    \exists \boldsymbol{x}\in\Omega
    \ \mathrm{such\ that}\
    N_i(\boldsymbol{x})N_j(\boldsymbol{x}) \neq 0 .
    \label{eq:node_graph}
\end{equation}

This definition means that pairs of nodes whose basis-function supports overlap are regarded as adjacent nodes.
The same definition can be used to construct nodal graphs for Lagrange basis functions, B-spline basis functions, and superimposed meshes containing both types of basis functions, as shown in \figurename~\ref{fig:node_graph}.

Equation~\eqref{eq:node_graph} defines pairs of nodes whose basis-function supports overlap as edges of the nodal graph.
In the coefficient matrices considered in this study, each matrix entry is generated from element integration associated with a corresponding pair of basis functions.
Therefore, only pairs of basis functions with overlapping nonzero supports can become candidates for nonzero matrix entries.
Consequently, the adjacency relation of the nodal graph corresponds to the candidate nonzero pattern of the coefficient matrices considered in this study.
Here, a candidate nonzero entry does not mean that the actual value of the entry is necessarily nonzero; rather, it denotes an entry that should be stored as a possible nonzero component in the sparse matrix.
In Section~\ref{sq:nodal_graph_par}, domain decomposition is performed based on this nodal graph.

\subsubsection{Construction of the nodal graph in BSFEM}

An important point is that the definition in Eq.~\eqref{eq:node_graph} is applicable regardless of whether nodes $v_i$ and $v_j$ belong to the same mesh.
Therefore, it can be applied not only to the classical finite element method, but also to BSFEM, which involves more complex mesh structures.

As described in Section~\ref{sq:BSFEM}, the BSFEM proposed by the authors uses a finite element mesh with B-spline basis functions and another finite element mesh with Lagrange basis functions.
The coefficient matrix of this method is decomposed into four submatrices:
$\boldsymbol{K}^{\mathrm{GG}}$, which accounts for products of B-spline basis functions;
$\boldsymbol{K}^{\mathrm{LL}}$, which accounts for products of Lagrange basis functions;
and $\boldsymbol{K}^{\mathrm{GL}}$ and $\boldsymbol{K}^{\mathrm{LG}}$, which account for products of B-spline and Lagrange basis functions.
Therefore, it is necessary to consider the adjacency relations between B-spline nodes, between Lagrange nodes, and between B-spline and Lagrange nodes.

In all these cases, as shown in \figurename~\ref{fig:node_graph}, two nodes $v_i$ and $v_j$ are regarded as adjacent when Eq.~\eqref{eq:node_graph} is satisfied.

\subsection{Domain decomposition of the nodal graph}
\label{sq:nodal_graph_par}

\begin{figure}
    \includegraphics[bb=0 0 1695 630,width=13cm]{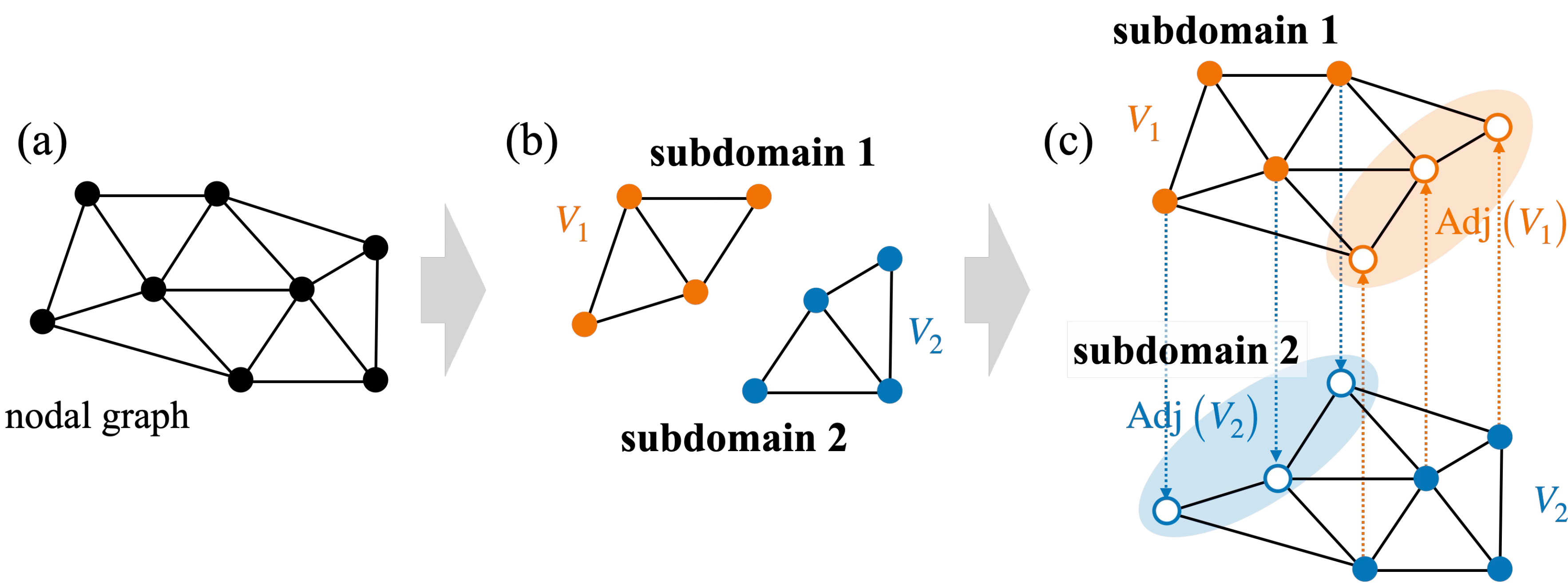}
    \centering
    \caption{Graph-based domain decomposition and definition of the adjacent node set.
    (a) Initial global nodal graph.
    (b) Partitioning of the graph into nonoverlapping subdomains $V_1$ and $V_2$ to balance computational weights.
    (c) Construction of extended subdomains.
    To compute the physical quantities on the internal nodes $V_k$ (solid circles), each subdomain must import data from adjacent nodes $\mathrm{Adj}(V_k)$ (hollow circles in shaded regions).
    This defines the required communication pattern without referring to the geometry.}
    \label{fig:node_graph_par}
\end{figure}

This subsection describes a domain decomposition method based on the nodal graph defined in the previous subsection.
A schematic illustration of a two-subdomain decomposition is shown in \figurename~\ref{fig:node_graph_par}.

\subsubsection{Load balancing using a weighted graph}
\label{sq:weighted_graph}

To achieve high scalability in finite element analysis on distributed-memory parallel computers, it is necessary to balance the computational load among processes while minimizing the amount of interprocess communication.
In graph theory, this problem is known as the minimum $k$-way cut problem \cite{karypis1998fast}.

This problem is formulated as the partitioning of the node set $V$ of a nodal graph $G^{\mathrm{node}}=(V,E)$ into $n_{\mathrm{dom}}$ subsets $V_k$ ($k=1,\dots,n_{\mathrm{dom}}$).
Each subset $V_k$ satisfies the following conditions:
\begin{align}
    V_k \cap V_l &= \emptyset \quad (k \neq l). \label{eq:direct_sum_dec1}
\end{align}

In the partitioning of a no-weight graph, the objective is to balance the number of nodes $|V_k|$ in each subdomain.
However, in the BSFEM considered in this study, the computational cost associated with each element differs significantly because of differences in the types of basis functions and integration algorithms used to generate the submatrices
$\boldsymbol{K}^{\mathrm{GG}}$, $\boldsymbol{K}^{\mathrm{LL}}$, $\boldsymbol{K}^{\mathrm{GL}}$, and $\boldsymbol{K}^{\mathrm{LG}}$.
Therefore, simply balancing the numbers of nodes or elements does not necessarily balance the computational load among processes.

To account for this nonuniformity in computational load, this study extends the nodal graph to a weighted graph.
Specifically, by appropriately distributing the element-wise computational cost to the associated nodes, a positive node weight $w_i^{\mathrm{V}}$, representing the computational cost, is assigned to each graph node $v_i$.
To keep the discussion focused, the element-wise computational cost used in this paper is detailed in Appendix~\ref{sec:appendix_A}, and the specific formulation of the node weight $w_i^{\mathrm{V}}$ is presented in Appendix~\ref{sec:appendix_B}.

The ideal domain decomposition problem is then redefined as an optimization problem that minimizes the total number of edges $J$ in the cut-edge set $E^{\mathrm{cut}}$, while balancing the total node weight $W_k^{\mathrm{V}}$ in each subdomain $k$:
\begin{align}
    &\text{minimize} \quad J = |E^{\mathrm{cut}}|, \label{eq:mincut1}\\
    &\text{subject to} \quad W^{\mathrm{V}}_k = W^{\mathrm{V}}_l \quad (\forall k, l). \nonumber
\end{align}
Here, the variables are defined as
\begin{align}
    W^{\mathrm{V}}_k &= \sum_{v_i \in V_k} w^{\mathrm{V}}_i, \\
    E^{\mathrm{cut}} &= \{ \left( v_i, v_j \right) \in E \mid \exists p, q \text{ such that } v_i \in V_p, v_j \in V_q, p \neq q \}.
\end{align}
Through this optimization, the internal node sets $V_k$ with appropriately distributed computational loads are determined, as shown in \figurename~\ref{fig:node_graph_par}(b). In the present implementation, each subdomain is assigned to one MPI process.

However, when the strict load-balance condition in Eq.~\eqref{eq:mincut1} is imposed on discrete node assignments, it may not be possible to exactly match the weights of all subdomains.
Moreover, for large-scale graphs, it is not practical to obtain an exact optimal partition that simultaneously minimizes the cut size and enforces strict load balance.
Therefore, in this study, the strict load-balance condition in Eq.~\eqref{eq:mincut1} is replaced by a relaxed constraint using an allowable imbalance ratio.
Specifically, we employ METIS \cite{karypis1998fast}, a representative graph partitioning library, and solve the relaxed problem given by Eq.~\eqref{eq:mincut2}:
\begin{align}
    &\text{minimize} \quad J = |E^{\mathrm{cut}}|, \label{eq:mincut2} \\
    &\text{subject to} \quad
    \frac{1}{1 + \epsilon} \frac{W_{\mathrm{total}}}{n_{\mathrm{dom}}}
    \le W^{\mathrm{V}}_k \le
    (1 + \epsilon) \frac{W_{\mathrm{total}}}{n_{\mathrm{dom}}}
    \quad (\forall k = 1, \dots, n_{\mathrm{dom}}). \nonumber
\end{align}
Here, $1+\epsilon$ is the allowable load imbalance ratio, and
$W_{\mathrm{total}} = \sum_{v_i \in V} w^{\mathrm{V}}_i$ is the total weight of the entire graph.
In all analyses performed in this study, the allowable load imbalance ratio was set to
$1+\epsilon=1.03 \ (\epsilon=0.03)$.

Although the actual implementation of the present method uses METIS, the proposed framework itself does not depend on a specific partitioning library. The essential idea of the proposed method is to represent both intra- and inter-mesh interactions among computational points in the superimposed B-spline and Lagrange meshes within a single nodal graph. This generalized graph representation enables existing graph partitioning methods to be applied to BSFEM without treating individual meshes separately.

\subsubsection{Subgraphs including adjacent nodes}

Each subset of nodes $V_k$ obtained by graph partitioning is a nonoverlapping set with respect to the other subdomains.
However, compared with the original single graph structure before partitioning, the edge information crossing subdomain boundaries is missing.
In finite element analysis, graph edges represent interactions, or data dependencies, between computational points.
Therefore, if such edges are cut, the information required for computation becomes incomplete, and the computation cannot be closed within each subdomain.

To address this issue, this study defines an extended region that complements the missing graph information so that each partitioned subdomain retains interaction information equivalent to that of the original single graph structure.
Specifically, for subdomain $k$, the set of nodes belonging to other subdomains $l$ ($l \neq k$) that are adjacent to the internal node set $V_k$, namely, the adjacent node set of $V_k$, is defined as
\begin{align}
    \mathrm{Adj}(V_k) = \{ v_j \in V \setminus V_k \mid \exists v_i \in V_k, (v_i, v_j) \in E \}.
\end{align}
A schematic illustration is shown in \figurename~\ref{fig:node_graph_par}(c).
In the upper-right panel of the figure, the adjacent node set $\mathrm{Adj}(V_1)$ for the internal node set $V_1$ of subdomain 1 is shown as open circles in the orange shaded region.
In the lower-right panel, the adjacent node set $\mathrm{Adj}(V_2)$ for the internal node set $V_2$ of subdomain 2 is shown as open circles in the blue shaded region.

Finally, the extended node set $\tilde{V}_k$ to be stored by subdomain $k$ is defined as the union of the internal node set and the adjacent node set:
\begin{align}
    \tilde{V}_k = V_k \cup \mathrm{Adj}(V_k).
    \label{eq:node_graph_k}
\end{align}
The subset of edges connecting nodes contained in the extended node set $\tilde{V}_k$ is denoted by $\tilde{E}_k$.
Each subdomain $k$ stores the extended subgraph $\tilde{G}^{\mathrm{node}}_k$ defined by these sets:
\begin{align}
    \tilde{G}^{\mathrm{node}}_k = (\tilde{V}_k, \tilde{E}_k).
\end{align}
This construction ensures that the data dependencies required for parallel computation, such as matrix generation and local matrix-vector products, are available within subdomain $k$.
As a result, computations in each subdomain can be performed without loss of information.

\subsubsection{Element assignment based on the partitioned nodal graph}

Element computations are assigned to each subdomain based on the extended node sets obtained from the partitioning of the nodal graph.
Let $V_e$ denote the element node set required to compute element $e$.
In this study, element $e$ is evaluated in subdomain $k$ when
\begin{equation}
    V_e \subseteq \tilde{V}_k
    \label{eq:elem_assignment}
\end{equation}
is satisfied.
That is, if all nodes required for the element computation are contained in the extended node set of subdomain $k$, the element is computed in that subdomain.

Elements located near subdomain boundaries may be evaluated redundantly in multiple subdomains.
However, contributions to the global matrix are assembled only into the row entries corresponding to the internal node set $V_k$ owned by each subdomain.
Therefore, consistency of the assembly result is preserved even when redundant element computations occur.

\subsubsection{Generality with respect to basis function degree}

\begin{figure}
    \includegraphics[bb=0 0 1373 469,width=13cm]{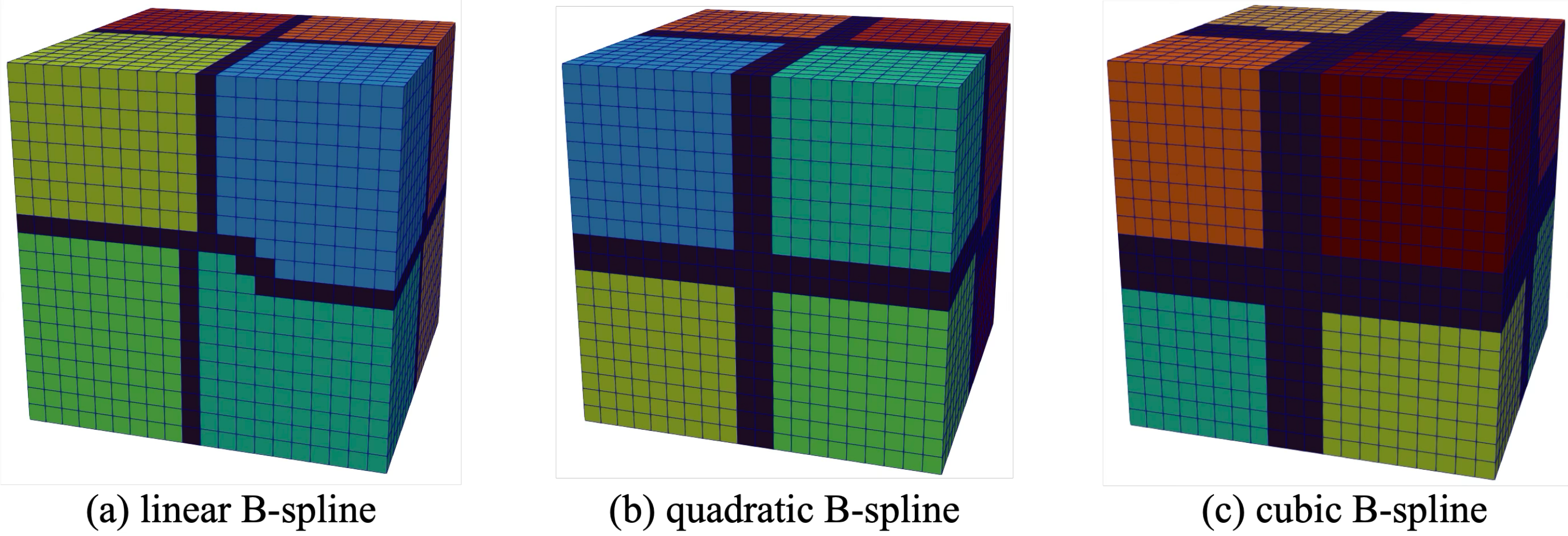}
    \centering
    \caption{Automatic adjustment of the number of shared-element layers based on B-spline degree. Domain decomposition results for (a) linear, (b) quadratic, and (c) cubic B-spline meshes partitioned into 8 domains. Black regions indicate shared elements between domains. The proposed graph-based method automatically generates the appropriate number of shared-element layers derived solely from the connectivity of the nodal graph.}
    \label{fig:b-spline_mesh_par}
\end{figure}

This subsection demonstrates that the proposed domain decomposition method based on the nodal graph can be applied through the same procedure regardless of the degree of the B-spline basis functions.
For B-spline basis functions, the support width increases as the degree becomes higher.
Consequently, the range of adjacent nodes required for matrix-vector products and element computations also expands according to the degree of the basis functions.

\figurename~\ref{fig:b-spline_mesh_par} shows the decomposition results for linear, quadratic, and cubic B-spline meshes with $20 \times 20 \times 20$ elements partitioned into 8 domains.
In the figure, elements included in multiple domains are shown in black.
Focusing on the number of element layers belonging to multiple domains near the subdomain boundaries, it can be observed that the number of shared-element layers automatically increases according to the basis function degree: one element layer for the linear case, two element layers for the quadratic case, and three element layers for the cubic case.

This behavior arises from the definition of the nodal graph in Eq.~\eqref{eq:node_graph}, which is based on the overlap of regions where the corresponding basis function values are nonzero.
A B-spline basis function of degree $p$ has compact support over $p+1$ knot spans, which correspond to $p+1$ elements when the multiplicity of each internal knot is one.
Therefore, as the degree $p$ increases, the basis function corresponding to a given node overlaps with basis functions associated with spatially more distant nodes, thereby expanding the adjacency range of that node in the nodal graph.
Because the graph partitioning algorithm determines the cut surfaces while accounting for this expanded adjacency relation, an appropriate number of shared-element layers corresponding to the degree $p$ is automatically constructed in physical space.

In the proposed method, domain decomposition can be performed for B-spline basis functions of arbitrary degree based on the same definition of the nodal graph.
That is, without introducing separate rules for constructing shared-element layers for individual basis functions, the data dependencies corresponding to the support width of the basis functions can be automatically handled as adjacency relations in the nodal graph.

\subsection{Linear algebra operations based on the nodal graph}

This subsection formulates parallel linear algebra operations based on the partitioned nodal graph defined in the previous subsection, with the aim of constructing a general parallel algorithm that is independent of the underlying numerical method.

As a solver for large-scale systems of linear equations, an iterative method represented by Krylov subspace methods is adopted from the viewpoint of parallel efficiency.
In this study, the matrix components assigned to each subdomain are stored in a block CSR format, and linear algebra operations are performed using this data structure.

The parallelization methods for the three principal elementary operations that constitute iterative solvers, namely vector addition, inner product, and matrix-vector multiplication, are described below on the extended nodal graph of subdomain $k$,
$\tilde{G}^{\mathrm{node}}_k = (\tilde{V}_k, \tilde{E}_k)$.

\subsubsection{Definition of the data structure}

In general physical simulations, each node has multiple physical quantities, such as the displacement components in the $x$, $y$, and $z$ directions in three-dimensional structural analysis.
Accordingly, the number of degrees of freedom per node is defined as $d$ ($d \ge 1$).

The vector $\boldsymbol{x}^k$ on subdomain $k$ is a block vector defined on the extended node set $\tilde{V}_k$.
The component $\boldsymbol{x}^k_i$ corresponding to each node $v_i$ is a $d$-dimensional vector:
\begin{align}
    \boldsymbol{x}^k_i \in \mathbb{R}^d \quad (v_i \in \tilde{V}_k).
\end{align}

The vector $\boldsymbol{x}^k$ over the entire subdomain is written in block form as the components $\boldsymbol{x}^k_{I}$ corresponding to the internal node set $V_k$ and the components $\boldsymbol{x}^k_{\Gamma}$ corresponding to the adjacent node set $\mathrm{Adj}(V_k)$:
\begin{align}
    \boldsymbol{x}^k = 
    \begin{pmatrix}
        \boldsymbol{x}^k_{I} \\
        \boldsymbol{x}^k_{\Gamma}
    \end{pmatrix}.
\end{align}

\subsubsection{Parallel computation of vector addition}

Consider the linear combination of two vectors $\boldsymbol{x}$ and $\boldsymbol{y}$, given by
$\boldsymbol{z} = \alpha \boldsymbol{x} + \beta \boldsymbol{y}$.
Each block component $\boldsymbol{z}^k_i$ of the local vector $\boldsymbol{z}^k$ in subdomain $k$ is computed as a linear operation between the corresponding $d$-dimensional vectors on each node:
\begin{align}
    \boldsymbol{z}^k_i
    =
    \alpha \boldsymbol{x}^k_i
    +
    \beta \boldsymbol{y}^k_i
    \quad (v_i \in V_k).
\end{align}
This operation is independent for each node and requires no interprocess communication.

\subsubsection{Parallel computation of vector inner products}

When computing the global inner product over the entire domain,
$\alpha = (\boldsymbol{x}, \boldsymbol{y})$,
a nonredundant summation must be performed.
From the definitions of domain decomposition in Eq.~\eqref{eq:direct_sum_dec1}, the complete node set $V$ is given by the disjoint union of the internal node sets $V_k$ of all subdomains:
\begin{align}
    V = \coprod_{k=1}^{n_{\mathrm{dom}}} V_k .
\end{align}
Therefore, the correct result is obtained by computing the local inner product only over the internal components, excluding the adjacent node components, and then aggregating the results over all subdomains:
\begin{align}
    \alpha = \sum_{k=1}^{n_{\mathrm{dom}}} \left( \sum_{v_i \in V_k} (\boldsymbol{x}^k_i)^{\top} \boldsymbol{y}^k_i \right).
\end{align}
In the present method, this operation is implemented using a global MPI communication routine, namely an allreduce operation.

\subsubsection{Parallel computation of matrix-vector products}

For the matrix-vector product $\boldsymbol{y} = \boldsymbol{K} \boldsymbol{x}$, subdomain $k$ is responsible for computing, or updating, only the components $\boldsymbol{y}^k_{I}$ corresponding to the internal node set $V_k$.
Each entry of the matrix $\boldsymbol{K}$ consists of a small dense matrix $\boldsymbol{K}_{ij}$ of size $d \times d$, where $d$ is the number of degrees of freedom per node, representing the action of node $v_j$ on node $v_i$.
Based on the adjacency relation of the nodal graph, this local operation can be written in the following block matrix form:
\begin{align}
    \boldsymbol{y}^k_{I} = \boldsymbol{K}^k_{II} \boldsymbol{x}^k_{I} + \boldsymbol{K}^k_{I \Gamma} \boldsymbol{x}^k_{\Gamma}.
    \label{eq:spmv_comp}
\end{align}
Here, $\boldsymbol{K}^k_{II}$ is a submatrix of size $d|V_k| \times d|V_k|$ that represents the interactions among nodes in the internal node set $V_k$, whereas $\boldsymbol{K}^k_{I \Gamma}$ is a submatrix of size $d|V_k| \times d|\mathrm{Adj}(V_k)|$ that represents the interactions between the internal node set $V_k$ and the adjacent node set $\mathrm{Adj}(V_k)$.

The first term in Eq.~\eqref{eq:spmv_comp},
$\boldsymbol{K}^k_{II}\boldsymbol{x}^k_I$,
can be computed using only data stored within the subdomain.
In contrast, the computation of the second term,
$\boldsymbol{K}^k_{I\Gamma}\boldsymbol{x}^k_{\Gamma}$,
requires the adjacent node components $\boldsymbol{x}^k_{\Gamma}$.
These components are copies of internal node components owned by neighboring subdomains.
Therefore, before performing the matrix-vector product, subdomain $k$ receives the corresponding components from the subdomains that own the nodes included in $\mathrm{Adj}(V_k)$, and updates $\boldsymbol{x}^k_{\Gamma}$ to the latest values.

\begin{figure}
    \includegraphics[bb=0 0 1712 801,width=11cm]{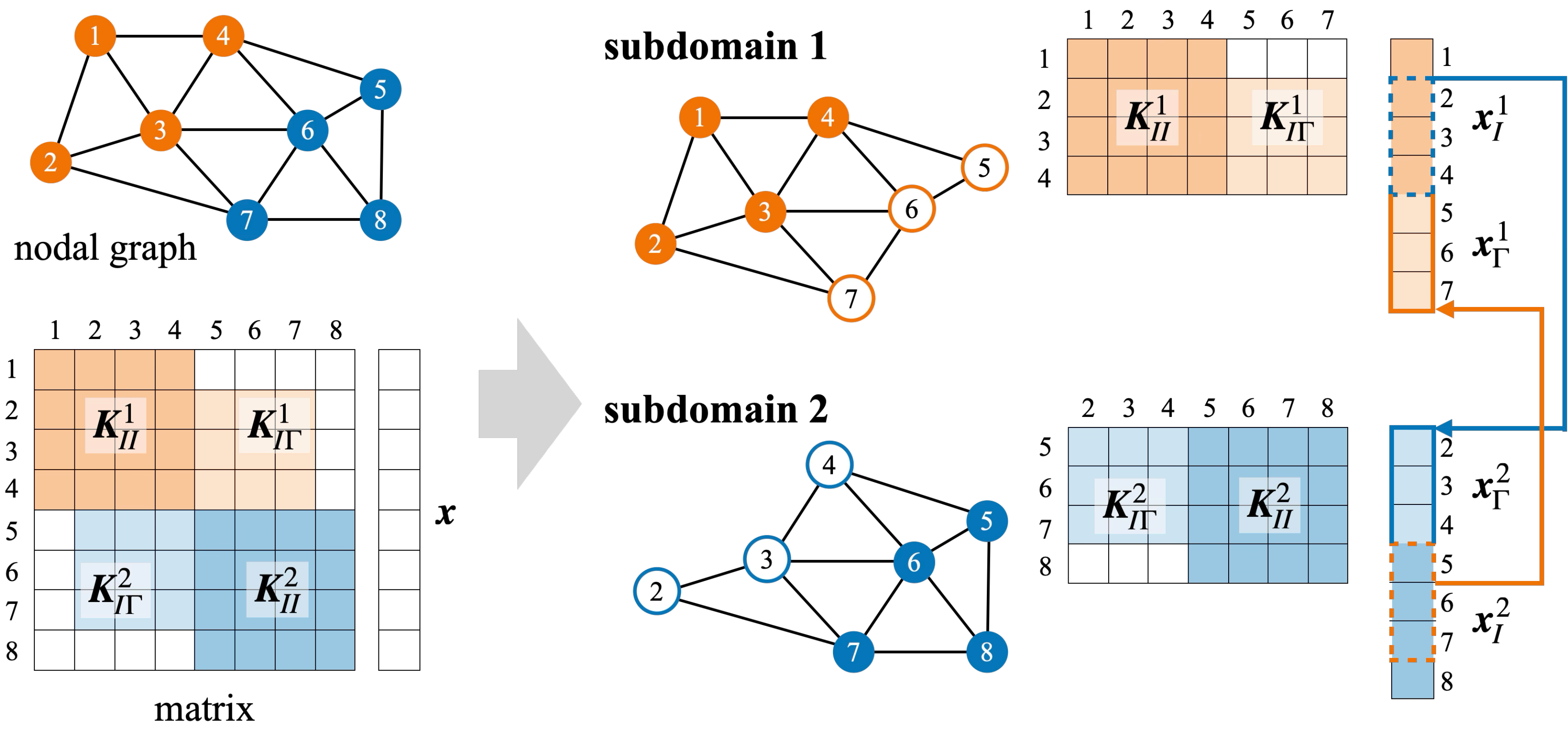}
    \centering
    \caption{Schematic representation of the parallel sparse matrix-vector product.
    The internal contribution $\boldsymbol{K}^k_{II}\boldsymbol{x}^k_I$
    is computed using only local data, whereas the adjacent contribution
    $\boldsymbol{K}^k_{I\Gamma}\boldsymbol{x}^k_{\Gamma}$
    requires communication with neighboring subdomains.}
    \label{fig:spmv_schematic}
\end{figure}

Accordingly, the matrix-vector product in subdomain $k$ is performed by first updating the adjacent node components $\boldsymbol{x}^k_{\Gamma}$ and then computing $\boldsymbol{y}^k_I$ according to Eq.~\eqref{eq:spmv_comp}.
A schematic illustration of this procedure is shown in \figurename~\ref{fig:spmv_schematic}.

\section{Verification of proposed method}
\label{sec:verification_parallelization}
In this section, the validity and parallel performance of the proposed graph-based domain-decomposition parallelization method are evaluated.
First, we verify that the proposed method does not degrade the accuracy of the numerical solution when the number of subdomains is varied.
This verification demonstrates that the ownership of degrees of freedom and elements and the MPI communication structure are appropriately constructed based on the proposed graph, thereby enabling domain-decomposition-based distributed-memory parallel computation of BSFEM.
Next, strong-scaling tests are performed for a fixed problem size to quantitatively evaluate the parallel performance and investigate the factors limiting the parallel performance.
Furthermore, as a demonstration of the utility of the proposed graph representation, cost-aware static load balancing is applied to the matrix-generation load imbalance identified through the parallel performance evaluation by incorporating computational costs into the node weights, and its effectiveness is examined.

\subsection{Verification conditions}
The Supercomputer for Quest to Unsolved Interdisciplinary Datascience (SQUID) system at Osaka University \citep{squid} was used for the present verification.
The hardware specifications of SQUID used for the parallel performance evaluation are summarized in \tablename~\ref{table:supacon}.

\begin{table}[th]
    \caption{Hardware specifications of the SQUID supercomputer system at Osaka University used for the parallel performance evaluation.}
    \label{table:supacon}
    \centering
    \small 
    \begin{tabular}{lr}
        \hline
        \rowcolor[rgb]{0.9, 0.9, 0.9}
        \multicolumn{2}{c}{Total system}\\
        Peak performance & 16.591 PFLOPS\\
        \hline
            \rowcolor[rgb]{0.9, 0.9, 0.9}
            \multicolumn{2}{c}{Compute node (general-purpose CPU nodes)}\\
            Number of nodes & 1,520 nodes \\
            Peak performance & 8.871 PFLOPS\\
            \multirow{2}{*}{CPU} & Intel${}^{\text{\textregistered}}$ Xeon${}^{\text{\textregistered}}$ Platinum 8368 \\
                                 & (Ice Lake / 2.4 GHz, 38 cores) $\times$ 2 \\
            Memory & 256 GB \\
        \hline
    \end{tabular}
\end{table}

The three-dimensional steady Poisson equation given in Eq.~\eqref{eq:poisson} is considered as the target problem.
The computational domain is defined as the three-dimensional cubic domain $\Omega^{\mathrm{G}}=[0,1]^3$, within which the local domain $\Omega^{\mathrm{L}}=[0,0.1]^3$ is superimposed.
The global and local meshes were constructed such that all elements were cubic.
Cubic B-spline basis functions were used for the global basis, whereas linear Lagrange basis functions were used for the local basis.
Based on the previous study \citep{magome2024higher}, the B-spline basis was constructed to maintain $C^2$ continuity across internal element boundaries.

Numerical integration in each element was performed using four-point Gauss quadrature in each coordinate direction.
The general-purpose parallel linear solver library Metagraph-Oriented Network for Linear Iterative Solvers (MONOLIS) \citep{monolis} was used to solve the system of linear equations.
The biconjugate gradient stabilized (BiCGSTAB) method was employed as the linear solver, with a relative residual tolerance of $1.0 \times 10^{-10}$.
The maximum number of iterations was set equal to the total number of degrees of freedom.
Diagonal scaling, which is independent of the number of subdomains, was employed as the preconditioner.

The error with respect to a known analytical solution was evaluated using the method of manufactured solutions \citep{roache1998verification}.
In this method, a sufficiently smooth function is prescribed as the analytical solution, and the corresponding source term is analytically determined such that the prescribed function satisfies the governing equation.
This allows the numerical solution to be directly compared with the known analytical solution and enables the consistency of the parallel implementation to be verified.

The manufactured solution used in this verification is given by
\begin{align}
    u = \sin(2 \pi x)\sin(2 \pi y)\sin(2 \pi z).
    \label{eq:manufactured_sol}
\end{align}
The corresponding Poisson equation is
\begin{align}
    \Delta u + 12 {\pi}^2 \sin(2 \pi x)\sin(2 \pi y)\sin(2 \pi z) = 0.
    \label{eq:manufactured_poisson}
\end{align}
As the boundary condition, the analytical values given by Eq.~\eqref{eq:manufactured_sol} were prescribed as Dirichlet conditions on the entire boundary of the computational domain.

The numerical error was evaluated using the relative error based on the $L^2$ norm, defined as
\begin{align} 
    \varepsilon_{L^2} 
    &= 
    \frac{
    \sqrt{
    \int_{\Omega^{\mathrm{G}}} 
    \left| \bar{u} - u \right|^2 d \Omega
    }
    }{
    \sqrt{
    \int_{\Omega^{\mathrm{G}}} 
    \left| u \right|^2 d \Omega
    }
    }
    =
    \frac{
    \sqrt{\varepsilon_{\mathrm{G}}+\varepsilon_{\mathrm{L}}}
    }{
    \sqrt{\varepsilon_{\mathrm{A}}}
    },
    \label{eq:sfem_error}\\
    \varepsilon_{\mathrm{G}}
    &=
    \int_{\Omega^{\mathrm{G}}\setminus\Omega^{\mathrm{L}}} 
    \left| \bar{u}^{\mathrm{G}} - u \right|^2 d \Omega, \nonumber \\
    \varepsilon_{\mathrm{L}}
    &=
    \int_{\Omega^{\mathrm{L}}} 
    \left| \bar{u}^{\mathrm{G}} + \bar{u}^{\mathrm{L}} - u \right|^2 d \Omega, \nonumber \\
    \varepsilon_{\mathrm{A}}
    &=
    \int_{\Omega^{\mathrm{G}}} 
    \left| u \right|^2 d \Omega. \nonumber
\end{align}
Here, $u$ denotes the manufactured solution, and $\bar{u}^{\mathrm{G}}$ and $\bar{u}^{\mathrm{L}}$ denote the numerical solutions obtained on the global and local meshes, respectively.
In addition, $\bar{u}$ denotes the SFEM solution and is defined as $\bar{u}=\bar{u}^{\mathrm{G}}$ in $\Omega^{\mathrm{G}}\setminus\Omega^{\mathrm{L}}$ and $\bar{u}=\bar{u}^{\mathrm{G}}+\bar{u}^{\mathrm{L}}$ in $\Omega^{\mathrm{L}}$.

In SFEM, domain integration is performed over multiple superimposed meshes; therefore, overlap detection between global and local elements is required.
In the present verification, the tolerance for overlap detection was set to $1.0 \times 10^{-8}$ to minimize the influence of errors associated with geometric detection.

\begin{table}[th]
    \centering
    \small
    \caption{Detailed configuration of the computational points and the number of elements for Mesh 1 used in the parallel performance evaluations.}
    \label{table:mesh1}
    \begin{tabular}{l *{6}{S[table-format=7.0, group-separator={,}]}} 
        \toprule
        & \multicolumn{3}{c}{Computational points} & \multicolumn{3}{c}{Elements} \\
        \cmidrule(lr){2-4} \cmidrule(lr){5-7}
        Mesh ID & {Total} & {Global} & {Local} & {Total} & {Global} & {Local} \\
        \midrule
        Mesh 1 & 166453 & 148877 & 17576 & 140625 & 125000 & 15625 \\
        \bottomrule
    \end{tabular}
\end{table}

\begin{figure}
    \centering
    \includegraphics[bb=0 0 845 551, width=8cm]{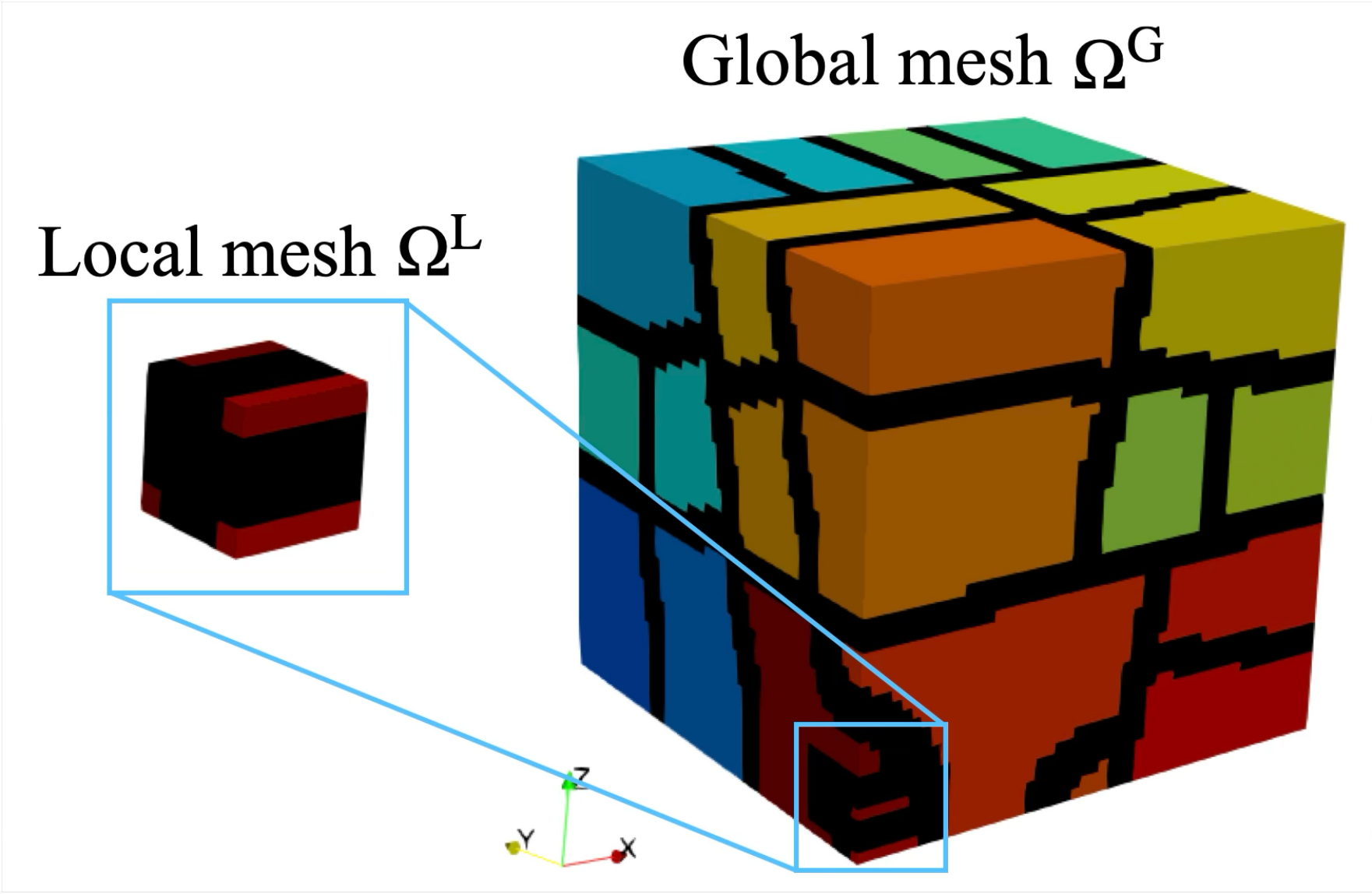}
    \caption{An example of a domain-decomposed mesh of Mesh 1 for BSFEM with 32 partitions. The black elements represent shared elements that span multiple subdomains. The node weights were uniformly set to 1 during domain decomposition.}
    \label{fig:sfem_ddm_mesh}
\end{figure}

\begin{figure}
    \centering
    \includegraphics[bb=0 0 410 310, width=10cm]{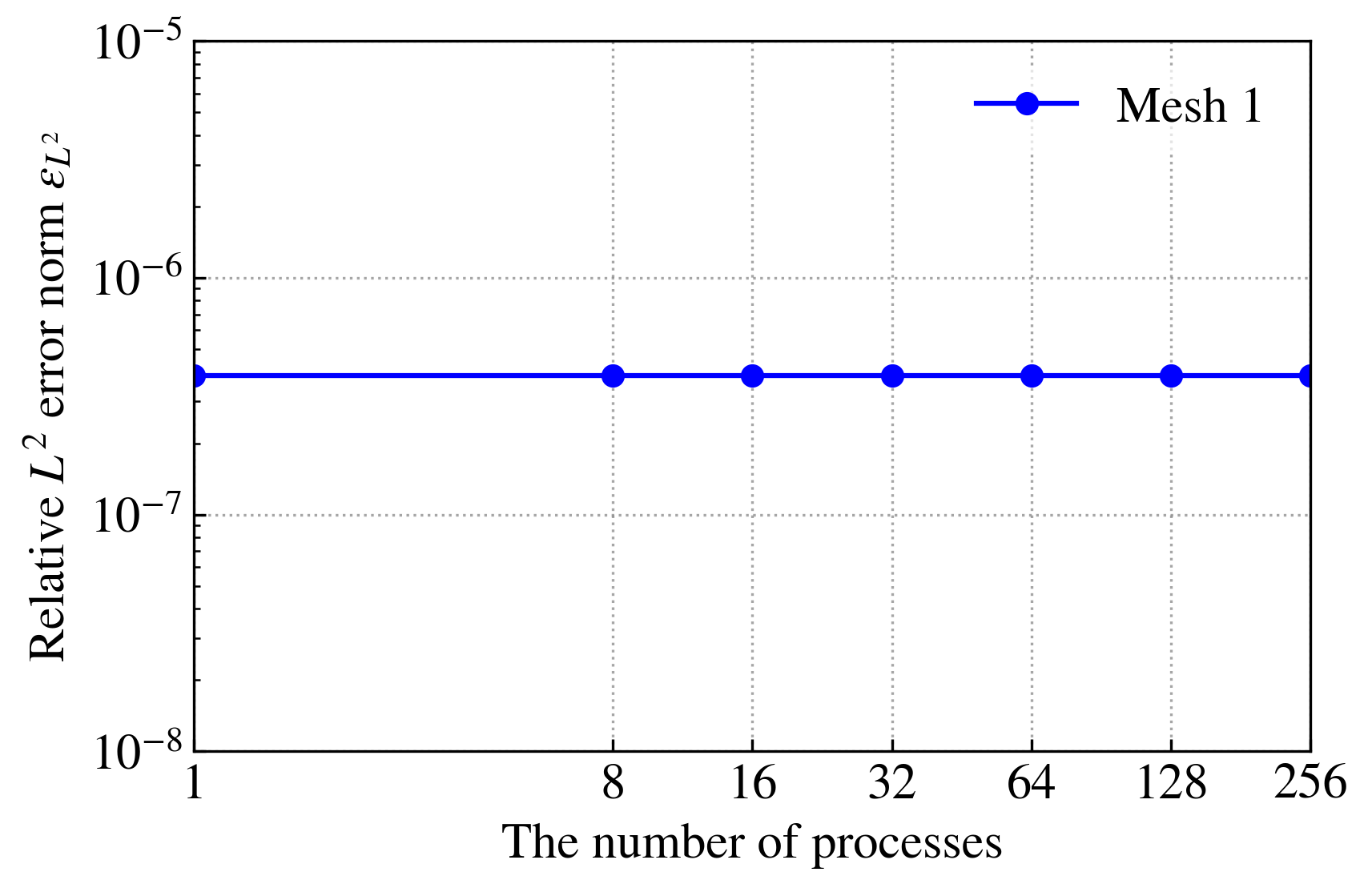}
    \caption{Relative $L^2$ error norm of the BSFEM solution as a function of the number of parallel processes. The results demonstrate that the solution accuracy is maintained as the number of subdomains is varied.}
    \label{fig:sfem_parallel_l2error}
\end{figure}

\subsection{Verification of graph-based domain-decomposition parallelization}

First, to confirm that the accuracy of the BSFEM solution is not degraded by the proposed generalized-graph-based domain-decomposition-based distributed-memory parallel computation, the solution accuracy was evaluated by varying the number of subdomains.

Mesh 1, shown in \tablename~\ref{table:mesh1}, was used in this evaluation.
The numbers of parallel processes were set to 1, 8, 16, 32, 64, 128, and 256, and the computations were performed on SQUID.
Considering the memory-bandwidth limitation and the memory capacity per compute node, the number of processes per compute node was limited to a maximum of 8.
Each parallel configuration was evaluated three times, and the average value was used as the result.
As the baseline condition in this subsection, the node weights of the nodal graph used for domain decomposition were uniformly set to 1.
For reference, an example of the domain decomposition with 32 processes is shown in \figurename~\ref{fig:sfem_ddm_mesh}.

The relative $L^2$ error defined in Eq.~\eqref{eq:sfem_error} was used as the evaluation metric.
The results are shown in \figurename~\ref{fig:sfem_parallel_l2error}.
The horizontal axis represents the number of parallel processes, and the vertical axis represents the relative $L^2$ error.
The relative errors are nearly identical for all numbers of processes, and no significant differences attributable to changes in the number of subdomains are observed.

Therefore, within the scope of this verification, the proposed domain-decomposition-based distributed-memory parallel computation was confirmed to preserve the accuracy of the BSFEM solution, with the solution accuracy remaining consistent as the number of subdomains was varied.

\subsection{Evaluation of parallel performance and computational load imbalance}

\begin{figure}
    \centering
    \includegraphics[bb=0 0 424 279, width=10cm]{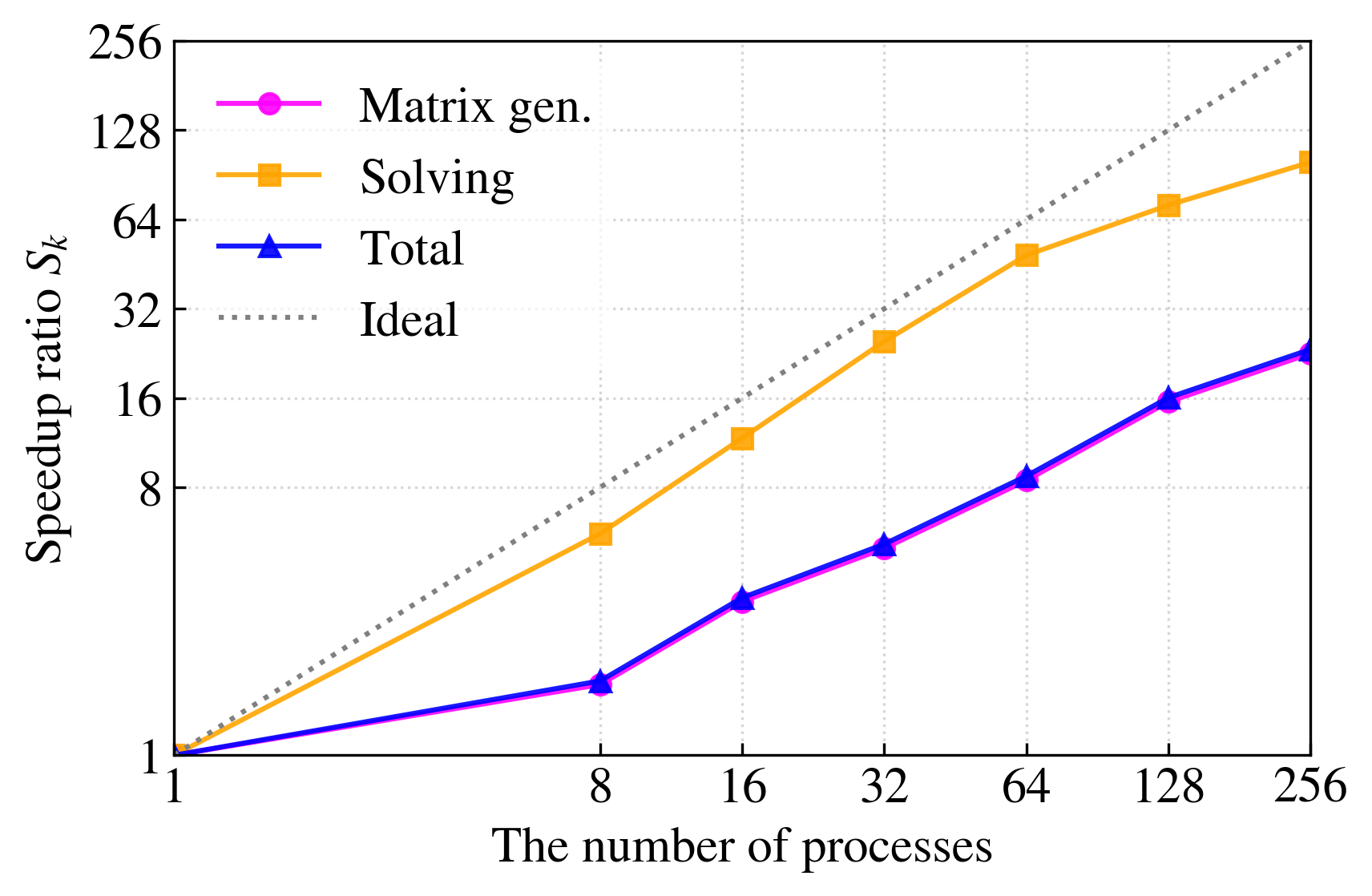}
    \caption{Strong-scaling performance: speedup $S_k$ for matrix generation time (magenta), linear solver time (orange), and total measured time (blue).}
    \label{fig:sfem_speedup}
\end{figure}

\begin{figure}
    \centering
    \includegraphics[bb=0 0 424 279, width=10cm]{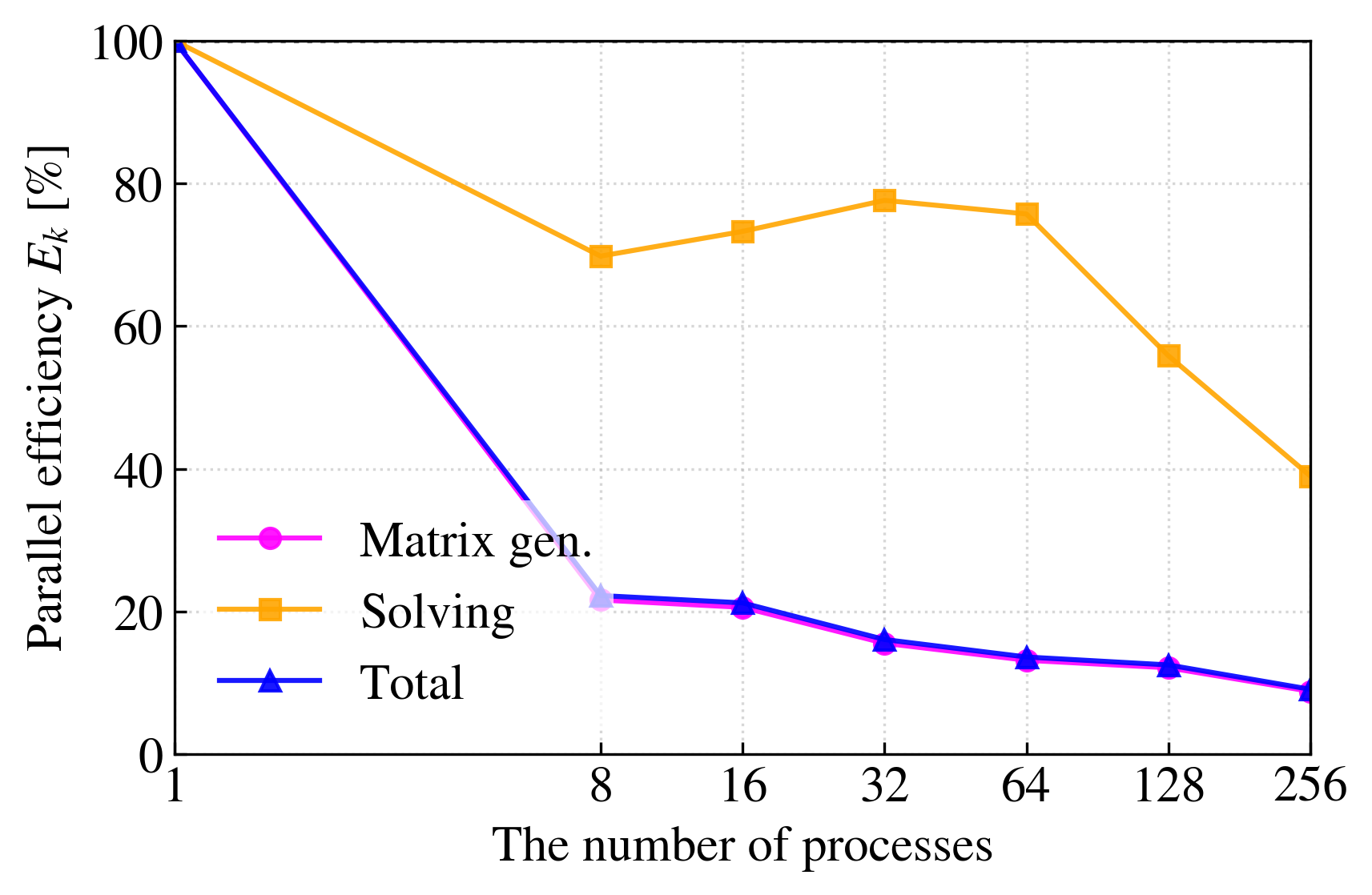}
    \caption{Strong-scaling performance: parallel efficiency $E_k$ for matrix generation time (magenta), linear solver time (orange), and total measured time (blue).}
    \label{fig:sfem_efficiency}
\end{figure}

\begin{table}[t]
    \centering
    \caption{Wall-clock times in the strong-scaling test for the no-weight case. $T_{\mathrm{Mat. gen.}}$ denotes the matrix generation time including synchronization waiting time, $T_{\mathrm{Solving}}$ denotes the linear solver time, and $T_{\mathrm{Total}}$ denotes the total measured time, defined as the sum of $T_{\mathrm{Mat. gen.}}$ and $T_{\mathrm{Solving}}$.}
    \label{table:strong_scaling_time}
    \begin{tabular}{
        S[table-format=3.0]
        |
        S[table-format=4.1]
        S[table-format=2.1]
        S[table-format=4.1]
        S[table-format=3.1]
    }
        \hline
        \multicolumn{1}{c|}{Number of processes $k$}
        & \multicolumn{1}{c}{$T_{\mathrm{Mat. gen.}}$ [s]}
        & \multicolumn{1}{c}{$T_{\mathrm{Solving}}$ [s]}
        & \multicolumn{1}{c}{$T_{\mathrm{Total}}$ [s]}
        & \multicolumn{1}{c}{$T_{\mathrm{Mat. gen.}}/T_{\mathrm{Solving}}$} \\
        \hline
        1   & 1302.9 & 55.4 & 1358.3 & 23.5 \\
        8   & 752.7  & 9.9  & 762.6  & 75.9 \\
        16  & 395.1  & 4.7  & 399.8  & 83.6 \\
        32  & 261.4  & 2.2  & 263.7  & 117.3 \\
        64  & 154.3  & 1.1  & 155.5  & 135.1 \\
        128 & 83.9   & 0.8  & 84.6   & 108.3 \\
        256 & 57.6   & 0.6  & 58.2   & 103.9 \\
        \hline
    \end{tabular}
\end{table}

Next, the parallel performance of the domain-decomposition-based distributed-memory parallel computation of BSFEM, whose numerical validity was confirmed in the previous subsection, is quantitatively evaluated through strong-scaling tests.
The factors limiting the parallel performance are then analyzed.

Mesh 1, as in the solution-accuracy evaluation in the previous subsection, was used in this verification.
The numbers of parallel processes, the number of measurements, the number of processes per compute node, and the node weights of the nodal graph used for domain decomposition were the same as those used in the previous evaluation.

The speedup $S_k$ defined in Eq.~\eqref{eq:s_k} and the parallel efficiency $E_k$ defined in Eq.~\eqref{eq:e_k} were used as the metrics for evaluating parallel performance.
Here, $k$ denotes the number of parallel processes, and $T_k$ denotes the measured time of the component under consideration using $k$ processes.
\begin{align}
    S_k &= \frac{T_1}{T_k},
    \label{eq:s_k}\\
    E_k &= \frac{S_k}{k}.
    \label{eq:e_k}
\end{align}

To evaluate the performance of the major computational components of the proposed domain-decomposed parallelization method, the time required for data input/output, basis-function construction in the reference space, and coupling-map search was excluded from the measurements.
The measured components were the matrix generation time (Matrix generation) and the time required to solve the system of linear equations (Solving).
Hereafter, the sum of these two components is referred to as the total measured time.

The calculated speedup $S_k$ is shown in \figurename~\ref{fig:sfem_speedup}, and the parallel efficiency $E_k$ is shown in \figurename~\ref{fig:sfem_efficiency}.
In both figures, the horizontal axis represents the number of parallel processes.
The magenta, orange, and blue lines represent the results for the matrix generation time, linear solver time, and total measured time, respectively.

These results show that, in the domain-decomposition-based distributed-memory parallel computation of BSFEM based on the proposed graph, both the matrix generation time and linear solver time decrease as the number of processes increases.
In particular, good scalability is obtained for solving the system of linear equations.
In contrast, the speedup for matrix generation is lower than that for the linear solver, thereby limiting the overall parallel performance.

The matrix generation time, linear solver time, and total measured time are listed in \tablename~\ref{table:strong_scaling_time}.
In the no-weight case, the linear solver time decreased from $T_{\mathrm{Solving}}=55.4$ s at $k=1$ to $0.6$ s at $k=256$, corresponding to a reduction of approximately $98.9\%$.
The matrix generation time also decreased from $T_{\mathrm{Mat. gen.}}=1302.9$ s at $k=1$ to $57.6$ s at $k=256$; however, the reduction was limited to approximately $95.6\%$.
Consequently, although matrix generation already accounted for a large fraction of the total measured time at low process counts, its relative contribution increased further at high process counts.
For example, at $k=64$, $T_{\mathrm{Mat. gen.}}/T_{\mathrm{Solving}}=135.1$, indicating that the matrix generation phase is the dominant factor limiting the overall parallel performance.

\begin{figure}
    \centering
    \begin{minipage}{0.49\textwidth}
        \centering
        \includegraphics[bb=0 0 1094 530, width=\textwidth]{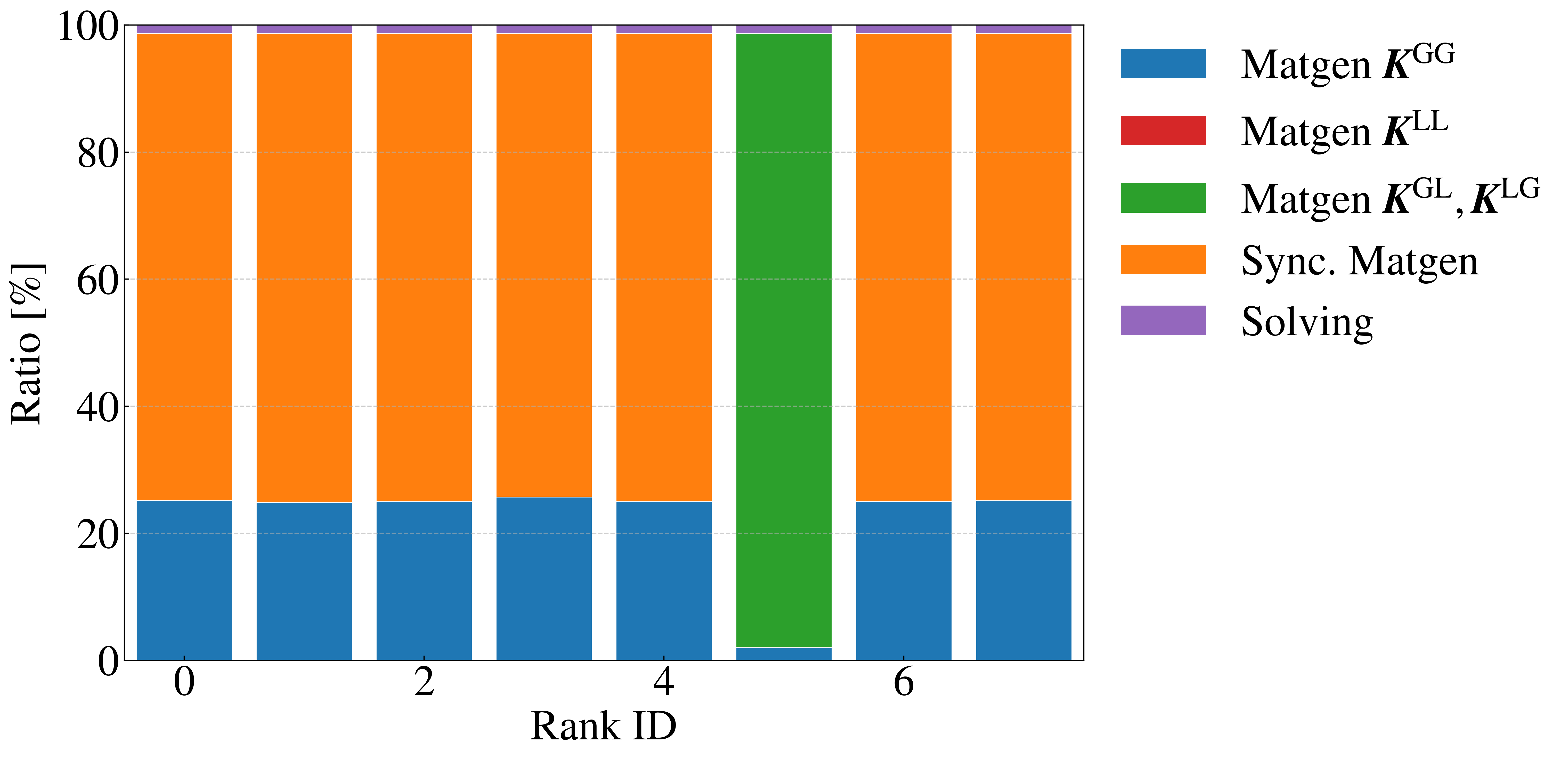}
        \par\small (a) 8 processes
    \end{minipage}%
    \hfill
    \begin{minipage}{0.49\textwidth}
        \centering
        \includegraphics[bb=0 0 1094 530, width=\textwidth]{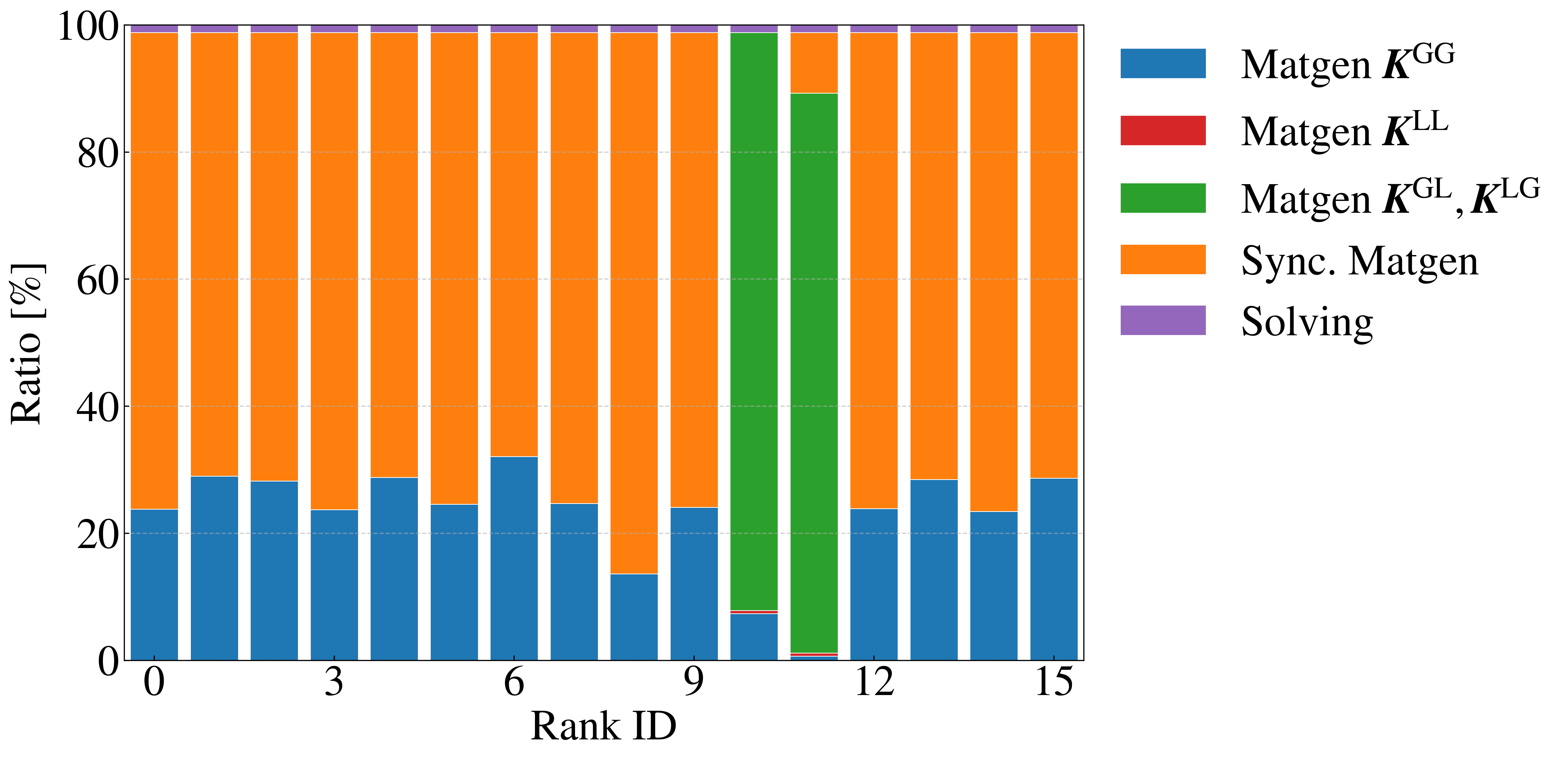}
        \par\small (b) 16 processes
    \end{minipage}

    \vspace{0.3cm}

    \begin{minipage}{0.49\textwidth}
        \centering
        \includegraphics[bb=0 0 1094 530, width=\textwidth]{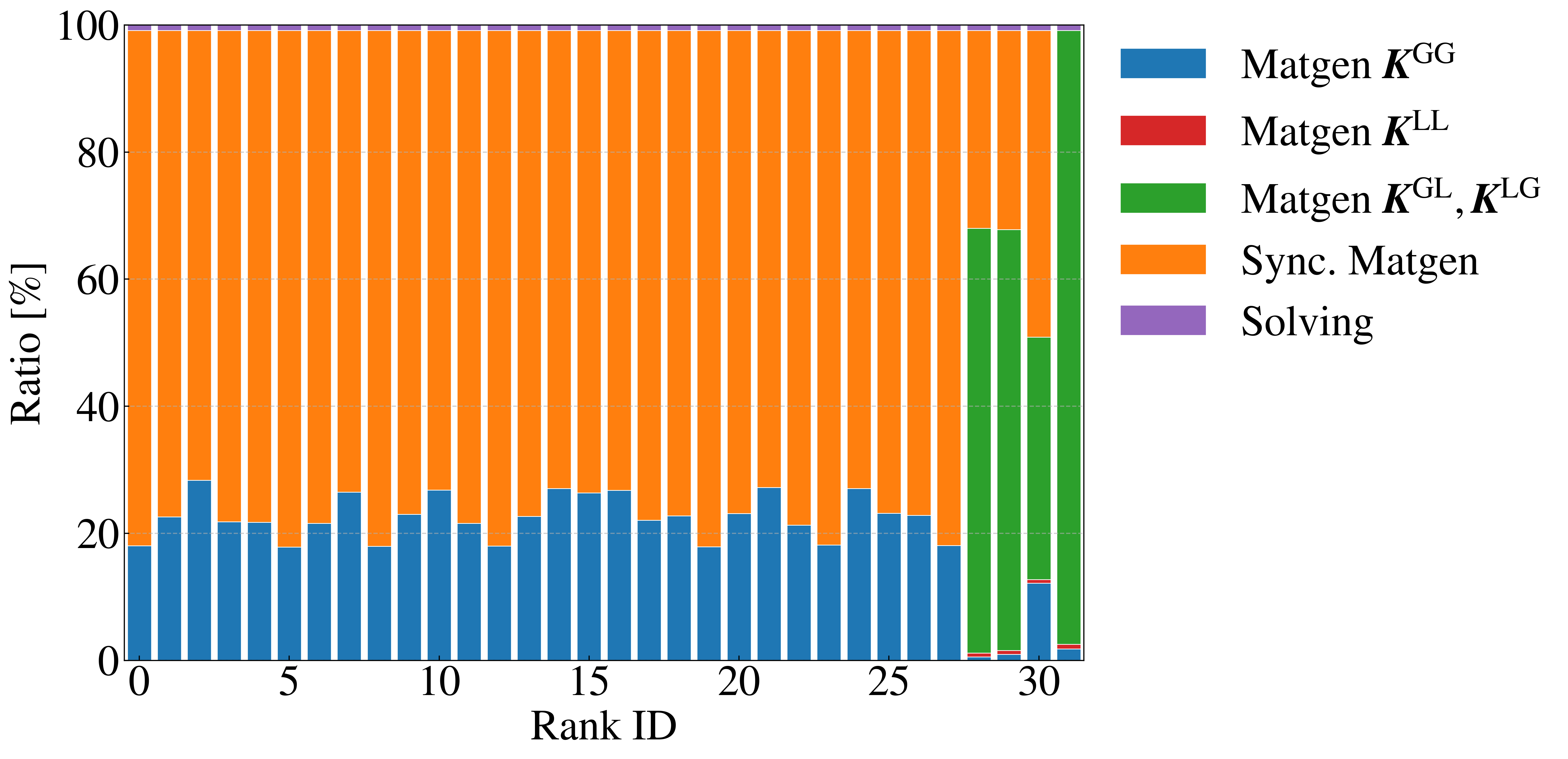}
        \par\small (c) 32 processes
    \end{minipage}%
    \hfill
    \begin{minipage}{0.49\textwidth}
        \centering
        \includegraphics[bb=0 0 1094 530, width=\textwidth]{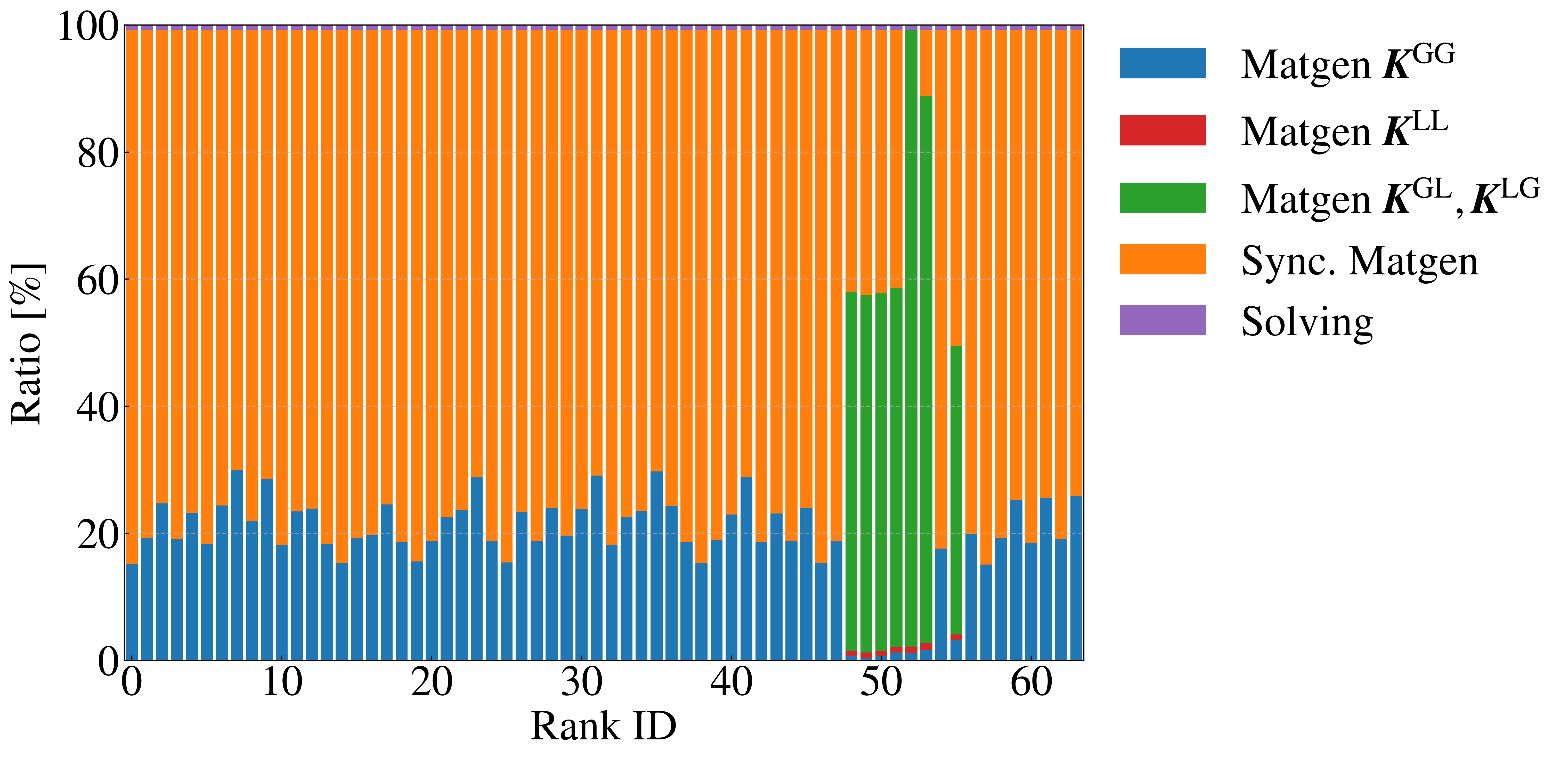}
        \par\small (d) 64 processes
    \end{minipage}

    \vspace{0.3cm}

    \begin{minipage}{0.49\textwidth}
        \centering
        \includegraphics[bb=0 0 1094 530, width=\textwidth]{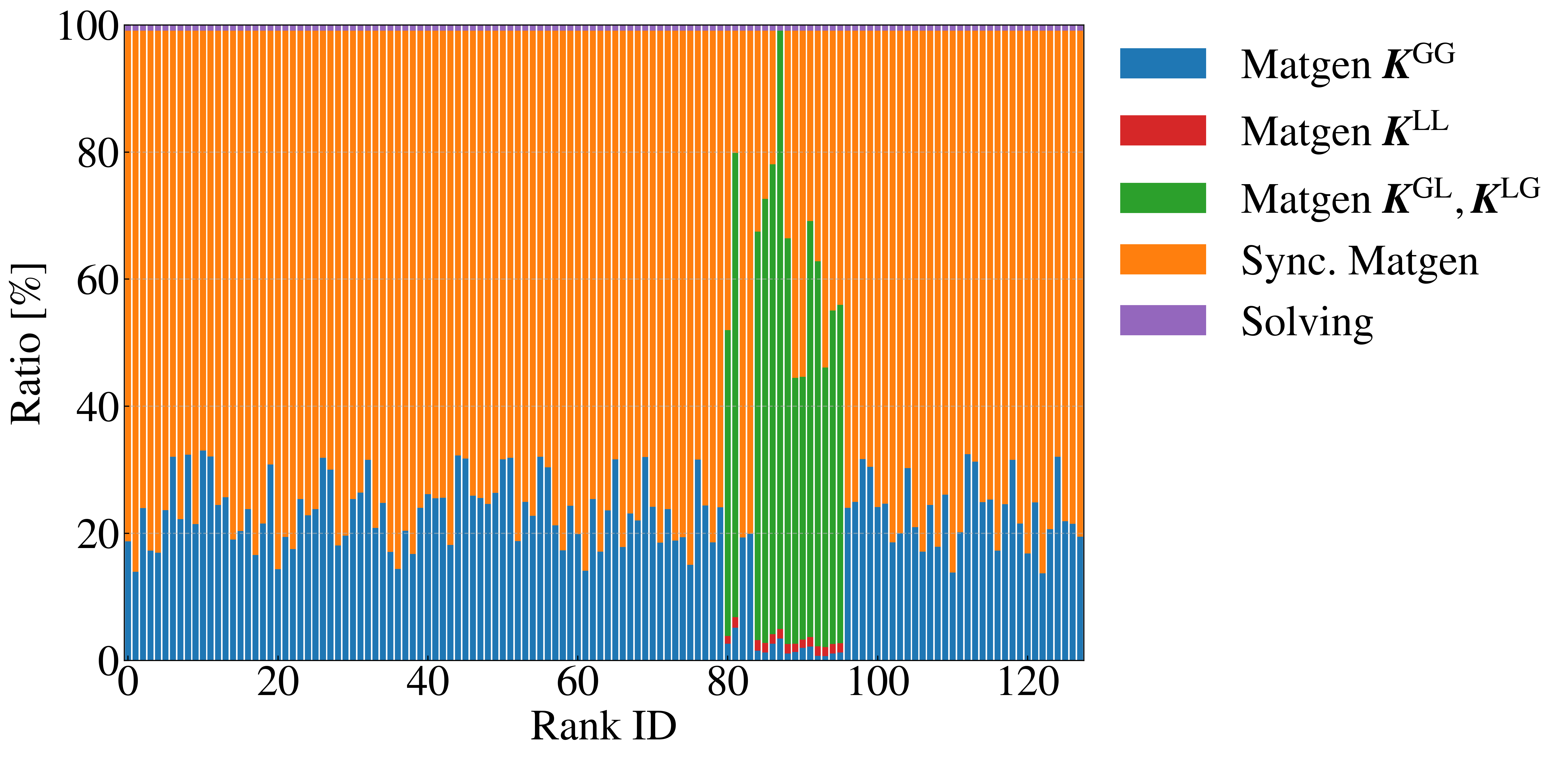}
        \par\small (e) 128 processes
    \end{minipage}%
    \hfill
    \begin{minipage}{0.49\textwidth}
        \centering
        \includegraphics[bb=0 0 1094 530, width=\textwidth]{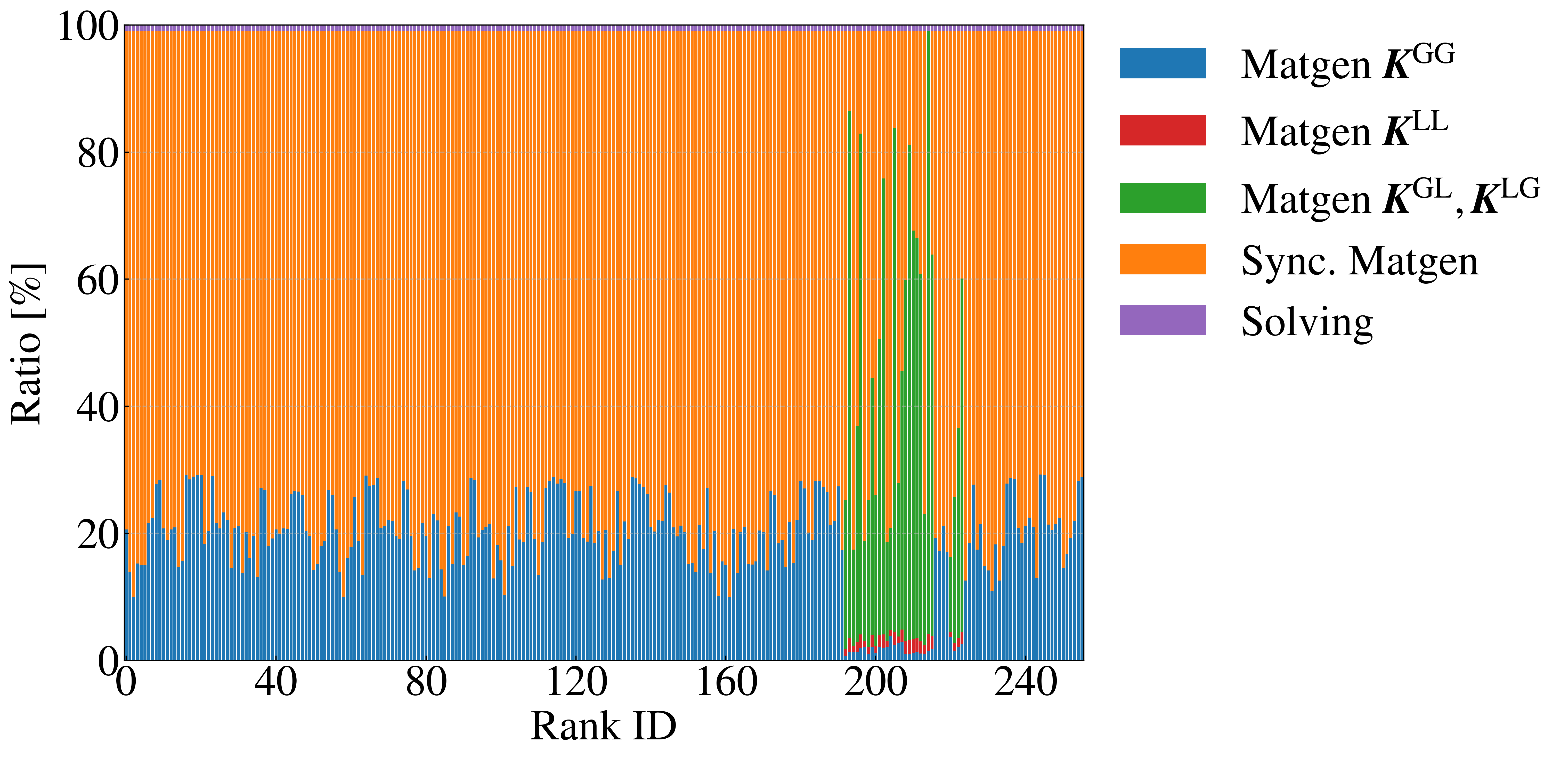}
        \par\small (f) 256 processes
    \end{minipage}

    \caption{Breakdown of the measured time for each process at different numbers of parallel processes. The blue, red, and green bars denote the generation times of the matrices $\boldsymbol{K}^{\mathrm{GG}}$, $\boldsymbol{K}^{\mathrm{LL}}$, and $\boldsymbol{K}^{\mathrm{GL}}/\boldsymbol{K}^{\mathrm{LG}}$, respectively. The orange bars denote the synchronization waiting time during the matrix generation phase, and the purple bars denote the linear solver time. The results indicate that the load imbalance in matrix generation, particularly in the generation of the coupling matrices, leads to substantial synchronization waiting time.}
    \label{fig:seven_time_ratio}
\end{figure}

To identify the cause of this bottleneck, the breakdown of the measured time for each process was investigated.
For the timing measurements, all processes were first synchronized using MPI\_Barrier before matrix generation.
The generation times of $\boldsymbol{K}^{\mathrm{GG}}$, $\boldsymbol{K}^{\mathrm{LL}}$, and $\boldsymbol{K}^{\mathrm{GL}}/\boldsymbol{K}^{\mathrm{LG}}$ were then measured separately for each process.
The generation time of each matrix was measured without intermediate synchronization, and all processes were synchronized again using MPI\_Barrier after all matrix generation processes had been completed.
In this study, the time from the completion of matrix generation by each process until that process passed through this synchronization point is defined as the synchronization waiting time in the matrix generation phase.

\figurename~\ref{fig:seven_time_ratio} shows the breakdown of the measured time for each process at each number of parallel processes.
The horizontal axis represents the process number, and the vertical axis represents the fraction of measured time.
The blue, red, and green bars represent the generation times of $\boldsymbol{K}^{\mathrm{GG}}$, $\boldsymbol{K}^{\mathrm{LL}}$, and $\boldsymbol{K}^{\mathrm{GL}}/\boldsymbol{K}^{\mathrm{LG}}$, respectively.
The orange bars represent the synchronization waiting time in the matrix generation phase, and the purple bars represent the linear solver time.

The figure shows that synchronization waiting time accounts for a substantial fraction of the matrix generation phase.
In particular, the generation of the coupling matrices $\boldsymbol{K}^{\mathrm{GL}}$ and $\boldsymbol{K}^{\mathrm{LG}}$ in BSFEM is concentrated in the local region $\Omega^{\mathrm{L}}$, where the global and local meshes overlap.
Consequently, when graph partitioning is performed using uniform node weights, this localized computational workload tends to be concentrated in a limited number of subdomains, resulting in synchronization waiting time for the other processes.
These results demonstrate that the matrix-generation workload associated with inter-mesh coupling in BSFEM is spatially nonuniform and that the resulting load imbalance is a major factor limiting parallel performance.

Therefore, to achieve practical parallel performance for BSFEM, the domain decomposition must account not only for the number of nodes but also for the computational workload associated with matrix generation.

\subsection{Demonstration of cost-aware static load balancing}

The previous subsection showed that, in domain-decomposition-based distributed-memory parallel computation of BSFEM, the matrix-generation workload associated with inter-mesh coupling is localized in the overlapping region and that domain decomposition using uniform node weights leads to load imbalance.
Because the proposed graph integrates the intra- and inter-mesh interactions in BSFEM into a single graph, BSFEM-specific matrix-generation costs can also be incorporated into the same partitioning representation.
In this subsection, as a demonstration of the utility of this graph representation, the matrix-generation costs of BSFEM are incorporated into the node weights to verify that static load balancing can account for the nonuniform computational workload.
The objective of this verification is not to propose a method for determining optimal node weights, but to demonstrate that load imbalance can be reduced by incorporating appropriate computational costs into the node weights.

In the node-weight model used in this study, the coefficients $\alpha_1$, $\alpha_2$, and $\alpha_3$ formulated in Appendix~\ref{sec:appendix_B} represent the computational workload associated with the generation of $\boldsymbol{K}^{\mathrm{GG}}$, $\boldsymbol{K}^{\mathrm{LL}}$, and $\boldsymbol{K}^{\mathrm{GL}}/\boldsymbol{K}^{\mathrm{LG}}$, respectively.

Among these components, the generation workloads of $\boldsymbol{K}^{\mathrm{GG}}$ and $\boldsymbol{K}^{\mathrm{LL}}$ are determined primarily by the number of elements, the degree of the basis functions, and the number of numerical integration points, as in element-matrix generation in conventional finite element methods.
In contrast, the results of the previous subsection indicated that the generation workload of the coupling matrices $\boldsymbol{K}^{\mathrm{GL}}$ and $\boldsymbol{K}^{\mathrm{LG}}$ is a major source of load imbalance.
However, the magnitude of this workload relative to the other matrix-generation workloads is not known a priori.
Therefore, in this verification, $\alpha_1$ and $\alpha_2$ were fixed, whereas only $\alpha_3$, which corresponds to the coupling-matrix generation workload, was varied to evaluate its effect on the workload distribution and the resulting changes in parallel performance.
Specifically,
\[
    \alpha_1 = 1, \quad
    \alpha_2 = 1, \quad
    \alpha_3 \in \{20, 40, 60, 80\}
\]
were used, and the workload distribution and parallel performance were compared for these four weighting conditions.
Increasing $\alpha_3$ corresponds to assigning a relatively larger computational workload to graph nodes involved in coupling-matrix generation during graph partitioning.
This is expected to distribute the high-workload region associated with the overlapping region over a larger number of subdomains and thereby reduce the synchronization waiting time in the matrix generation phase.
The total measured time was used as the evaluation metric.
For $k=1$, graph partitioning has no effect; therefore, the no-weight value was used as the reference for all weighting conditions.

\begin{figure}
    \centering
    \begin{minipage}{0.49\textwidth}
        \centering
        \includegraphics[bb=0 0 790 390, width=\textwidth]{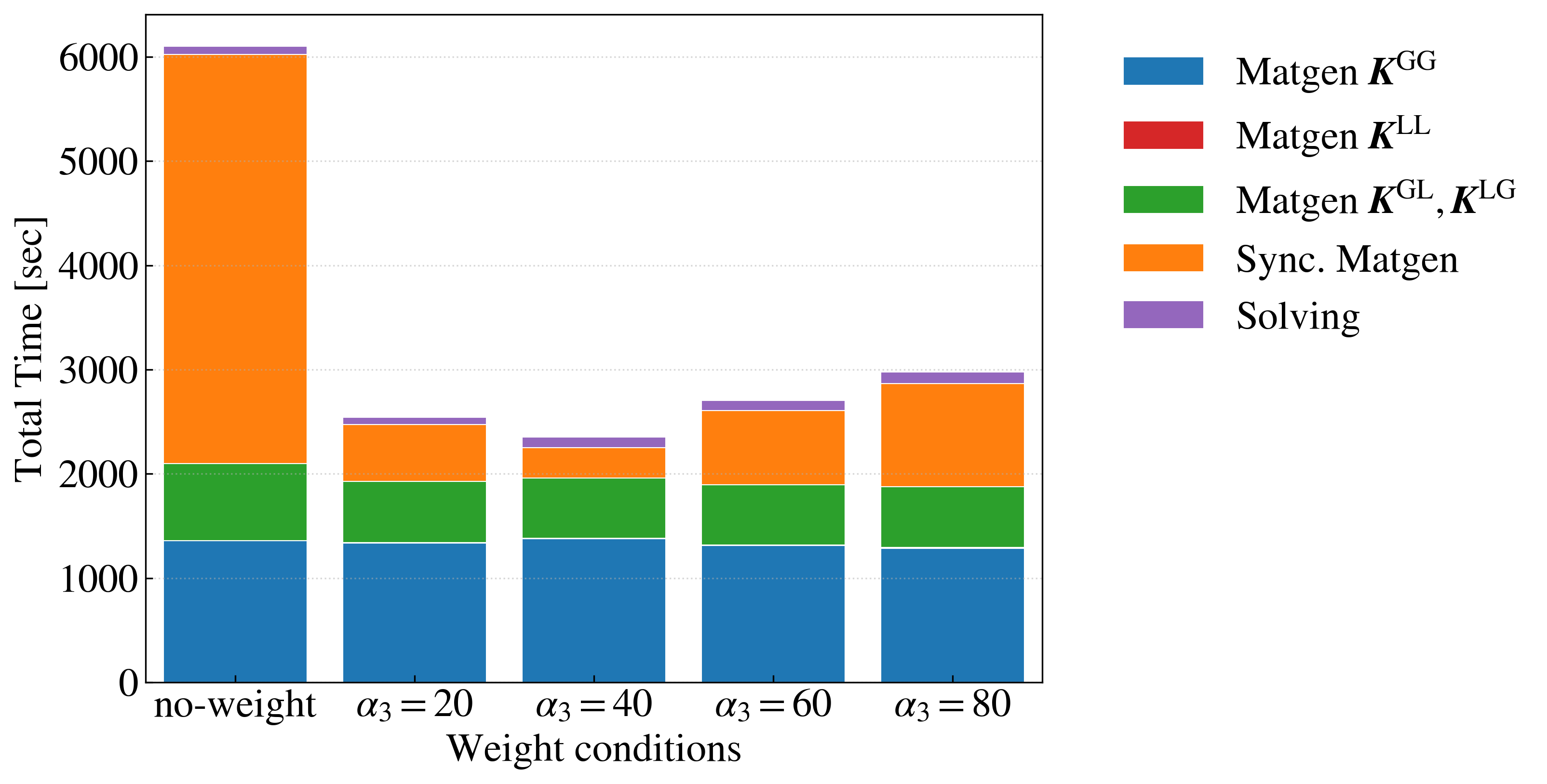}
        \par\small (a) 8 processes
    \end{minipage}%
    \hfill
    \begin{minipage}{0.49\textwidth}
        \centering
        \includegraphics[bb=0 0 790 390, width=\textwidth]{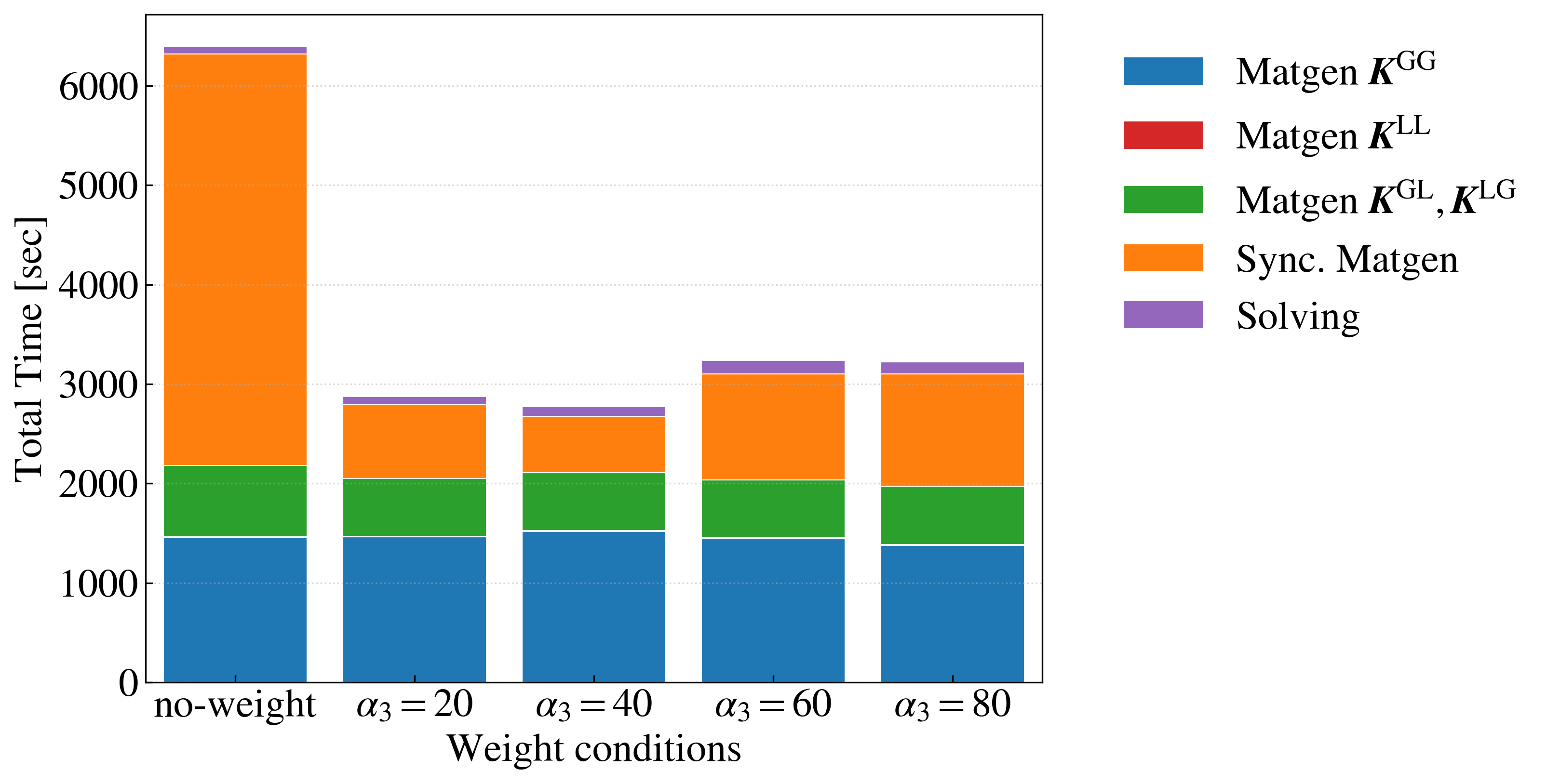}
        \par\small (b) 16 processes
    \end{minipage}

    \vspace{0.3cm}

    \begin{minipage}{0.49\textwidth}
        \centering
        \includegraphics[bb=0 0 790 390, width=\textwidth]{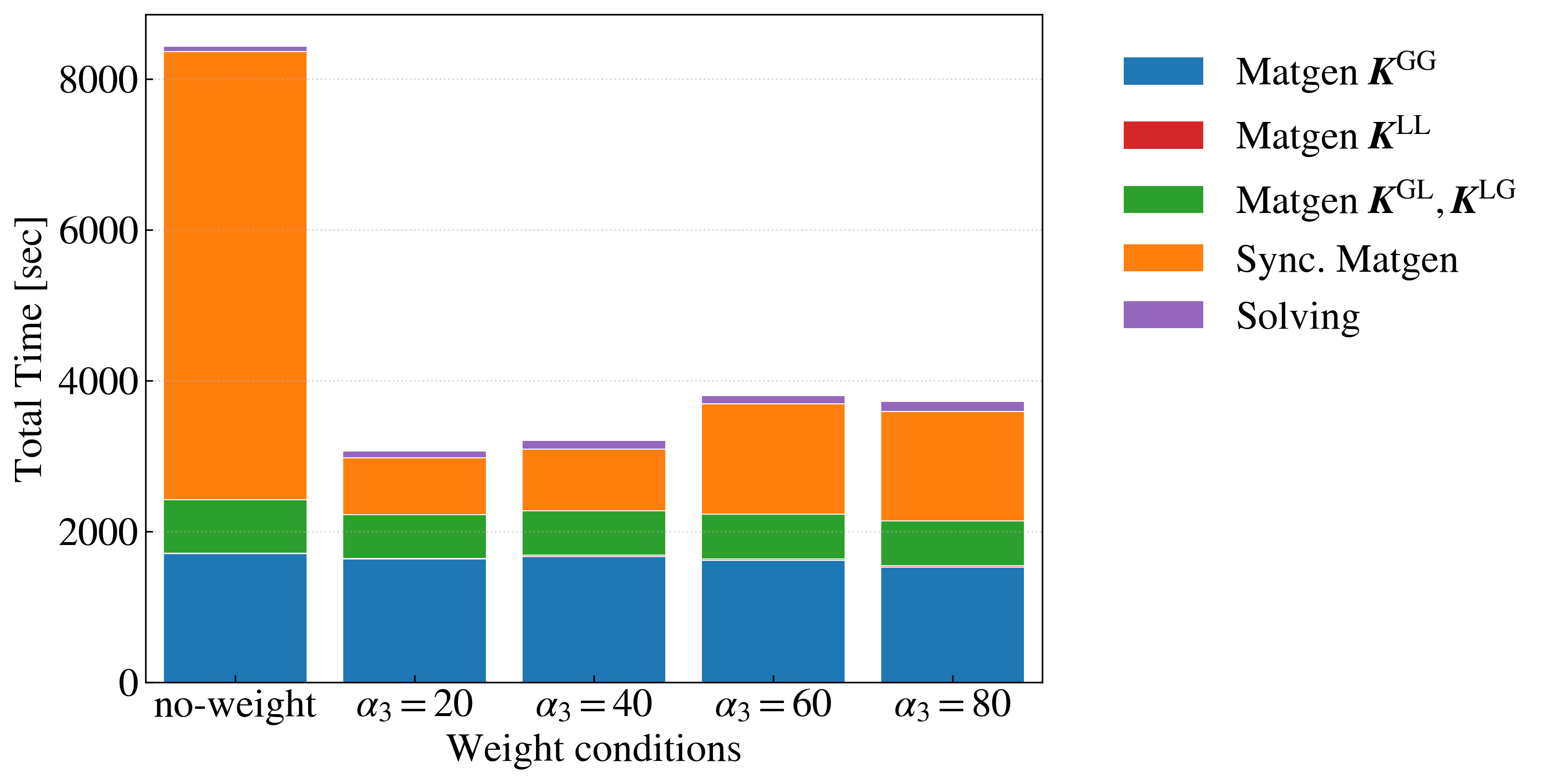}
        \par\small (c) 32 processes
    \end{minipage}%
    \hfill
    \begin{minipage}{0.49\textwidth}
        \centering
        \includegraphics[bb=0 0 790 390, width=\textwidth]{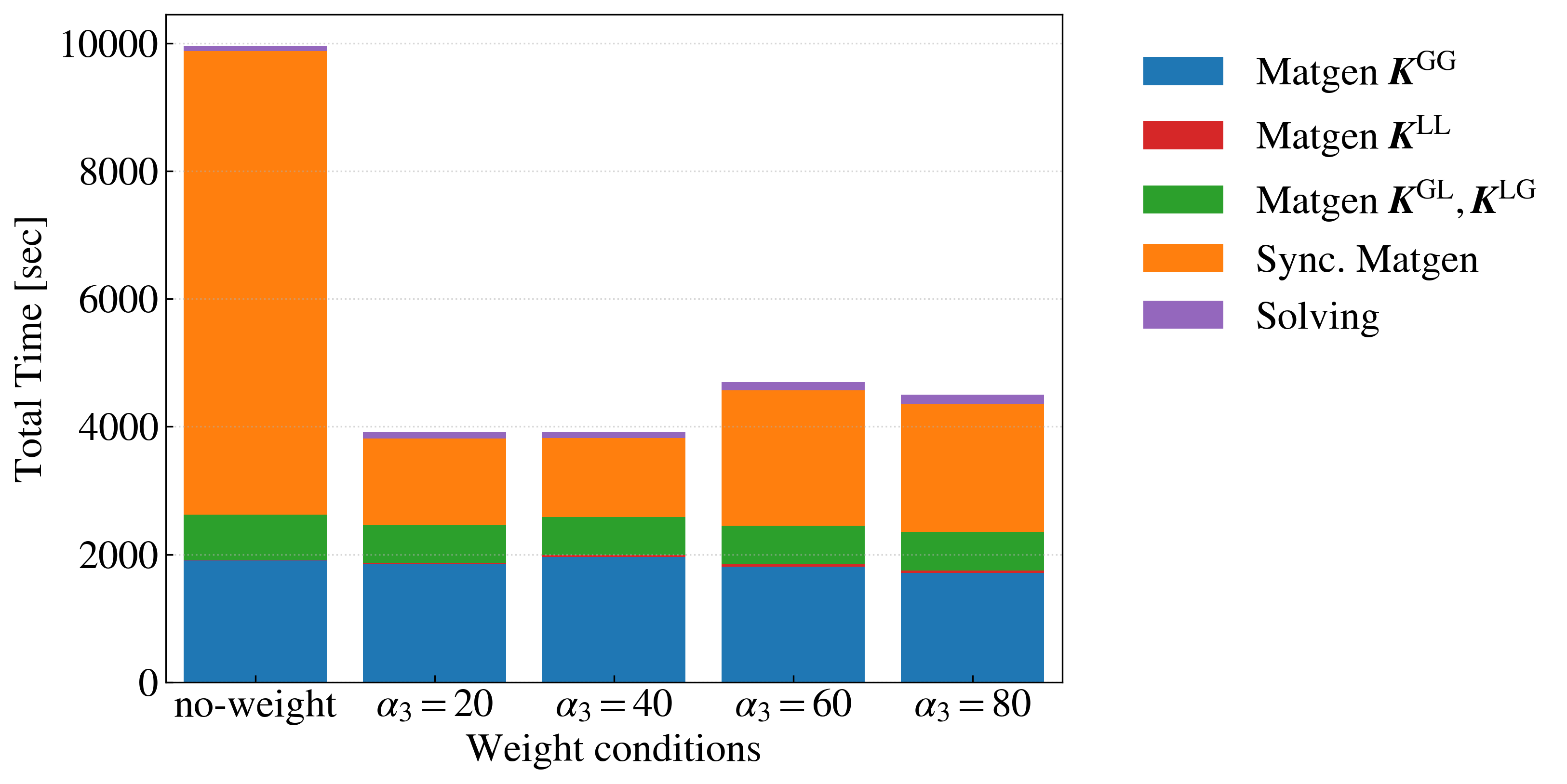}
        \par\small (d) 64 processes
    \end{minipage}

    \vspace{0.3cm}

    \begin{minipage}{0.49\textwidth}
        \centering
        \includegraphics[bb=0 0 790 390, width=\textwidth]{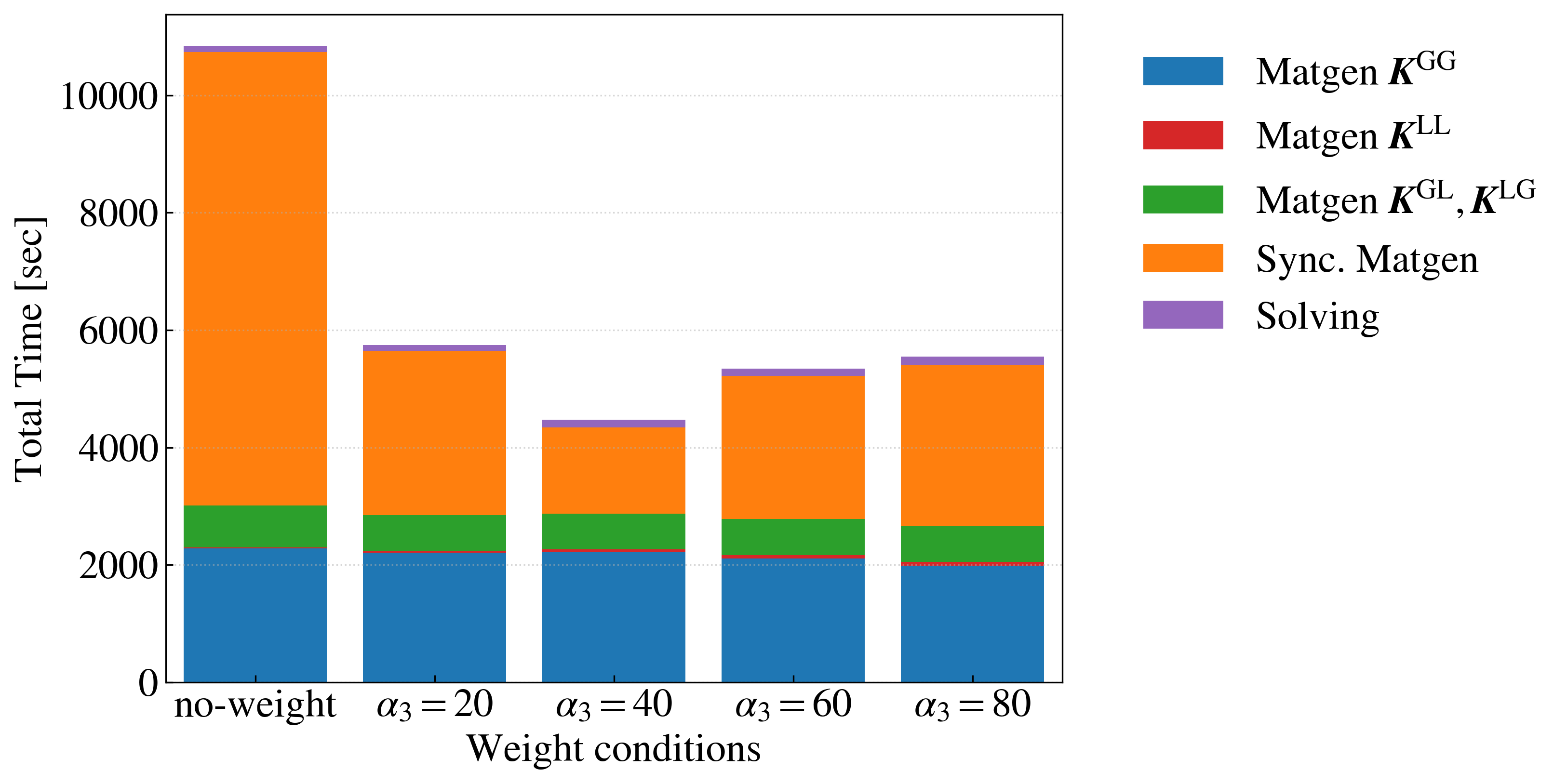}
        \par\small (e) 128 processes
    \end{minipage}
    \hfill
    \begin{minipage}{0.49\textwidth}
        \centering
        \includegraphics[bb=0 0 790 390, width=\textwidth]{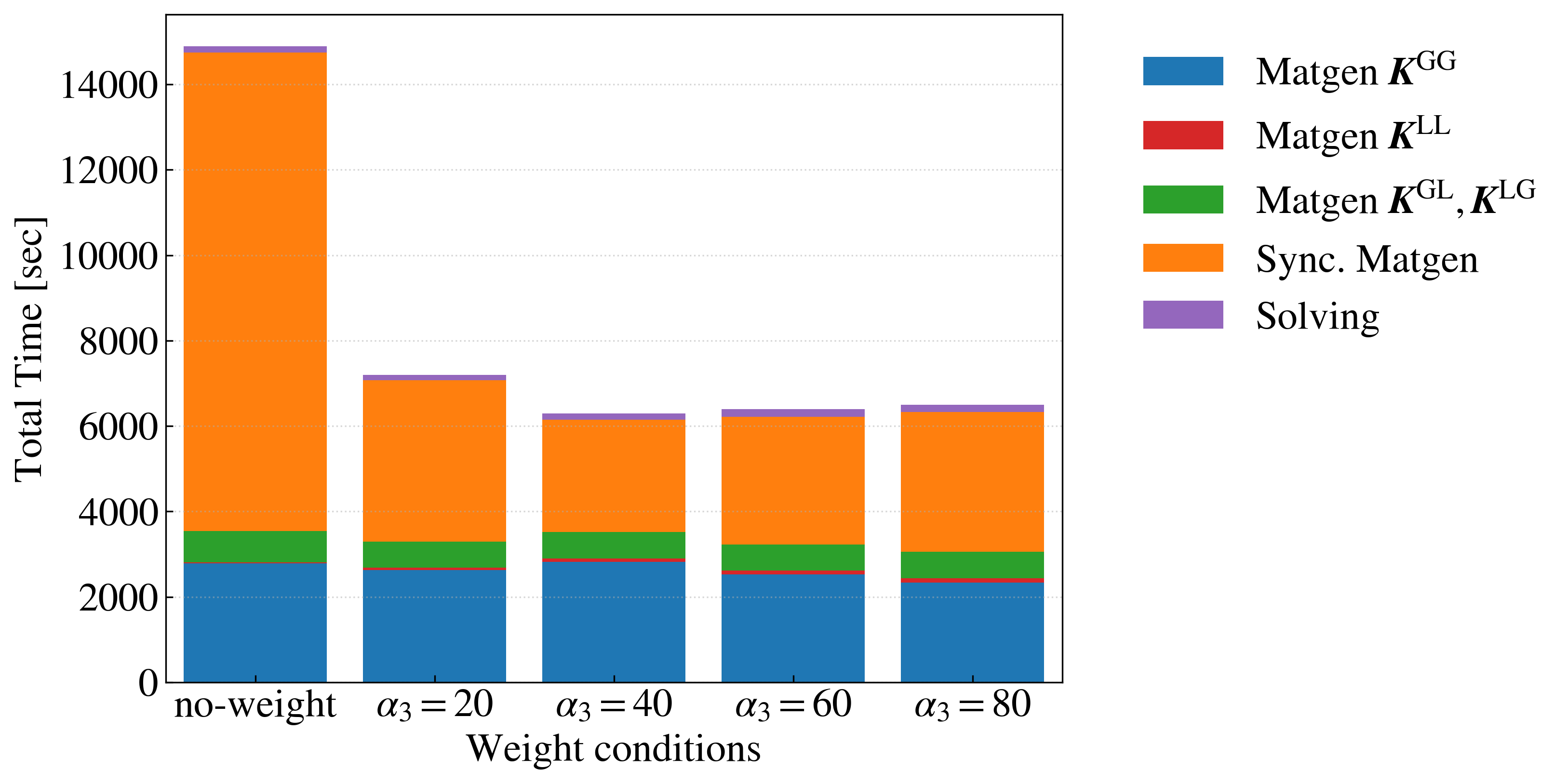}
        \par\small (f) 256 processes
    \end{minipage}

    \vspace{0.3cm}

    \caption{Breakdown of the accumulated measured time over all processes under different values of the weighting parameter $\alpha_3$. The blue, red, and green components denote the generation times of $\boldsymbol{K}^{\mathrm{GG}}$, $\boldsymbol{K}^{\mathrm{LL}}$, and $\boldsymbol{K}^{\mathrm{GL}}/\boldsymbol{K}^{\mathrm{LG}}$, respectively. The orange component denotes the synchronization waiting time during the matrix generation phase, and the purple component denotes the linear solver time. The results show that incorporating the coupling-matrix generation cost into the node weights substantially reduces the synchronization waiting time.}
    \label{fig:stacked_bar_alpha3}
\end{figure}

\begin{table}[t]
    \centering
    \caption{Total measured time under different weighting parameters.}
    \label{table:weighting_total_time}
    \begin{tabular}{
        S[table-format=3.0]
        |
        S[table-format=3.1]
        S[table-format=3.1]
        S[table-format=3.1]
        S[table-format=3.1]
        S[table-format=3.1]
    }
        \hline
        \multicolumn{1}{c|}{Number of processes $k$}
        & \multicolumn{1}{c}{No-weight}
        & \multicolumn{1}{c}{$\alpha_3=20$}
        & \multicolumn{1}{c}{$\alpha_3=40$}
        & \multicolumn{1}{c}{$\alpha_3=60$}
        & \multicolumn{1}{c}{$\alpha_3=80$} \\
        \hline
        8   & 762.6 & 317.8 & 293.9 & 337.8 & 372.1 \\
        16  & 399.8 & 179.6 & 173.1 & 202.3 & 201.4 \\
        32  & 263.7 & 95.9  & 100.2 & 118.9 & 116.4 \\
        64  & 155.5 & 61.1  & 61.3  & 73.4  & 70.3  \\
        128 & 84.6  & 44.9  & 34.9  & 41.8  & 43.3  \\
        256 & 58.2  & 28.1  & 24.6  & 25.0  & 25.4  \\
        \hline
    \end{tabular}
\end{table}

First, the effect of introducing node weights on the breakdown of measured time was evaluated.
For each number of parallel processes, the accumulated measured time over all processes was calculated, and its breakdown was compared.
The measured time components considered here consist of the generation times of the individual matrices, the synchronization waiting time in the matrix generation phase, and the linear solver time.
The results are shown in \figurename~\ref{fig:stacked_bar_alpha3}.
The horizontal axis represents the weighting condition, and the vertical axis represents the accumulated measured time over all processes.

The blue, red, and green components represent the generation times of $\boldsymbol{K}^{\mathrm{GG}}$, $\boldsymbol{K}^{\mathrm{LL}}$, and $\boldsymbol{K}^{\mathrm{GL}}/\boldsymbol{K}^{\mathrm{LG}}$, respectively.
The orange component represents the synchronization waiting time in the matrix generation phase, whereas the purple component represents the linear solver time.

As shown in \figurename~\ref{fig:stacked_bar_alpha3}, incorporating the matrix-generation costs into the node weights substantially reduces the synchronization waiting time between processes without substantially changing the accumulated matrix-generation time over all processes.
This indicates that the proposed weighting improves static load balancing.

The reduction in synchronization waiting time also results in a decrease in the total measured time.
As shown in \tablename~\ref{table:weighting_total_time}, for example, at $k=32$, the total measured time was $263.7$ s in the no-weight case, whereas it decreased to $100.2$ s for $\alpha_3=40$.
This corresponds to a reduction of $163.5$ s, or approximately $62.0\%$.

\begin{figure}
    \centering
    \begin{minipage}{0.49\textwidth}
        \centering
        \includegraphics[bb=0 0 1094 530, width=\textwidth]{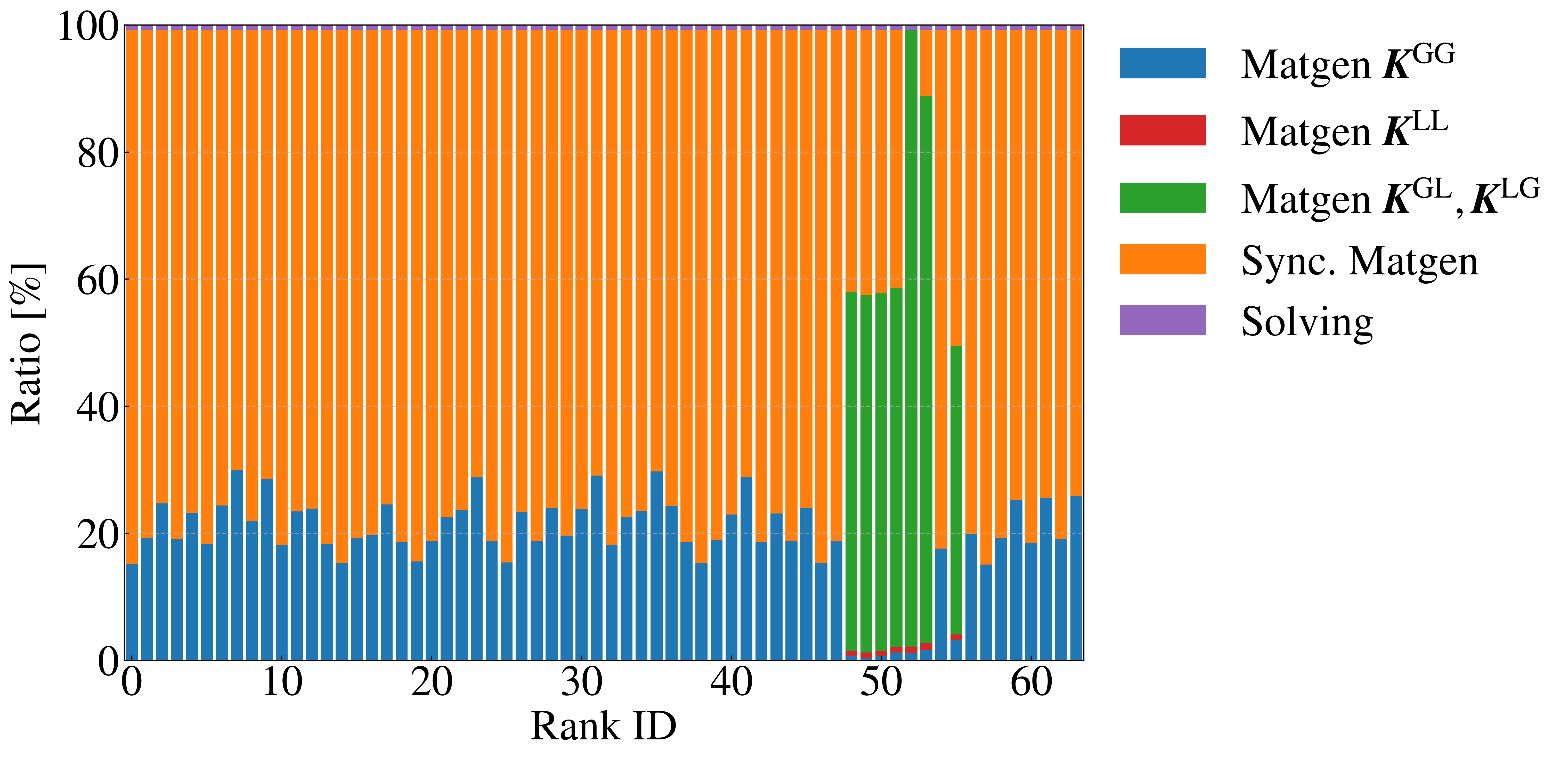}
        \par\small (a) No-weight
    \end{minipage}%

    \vspace{0.3cm}

    \begin{minipage}{0.49\textwidth}
        \centering
        \includegraphics[bb=0 0 1094 530, width=\textwidth]{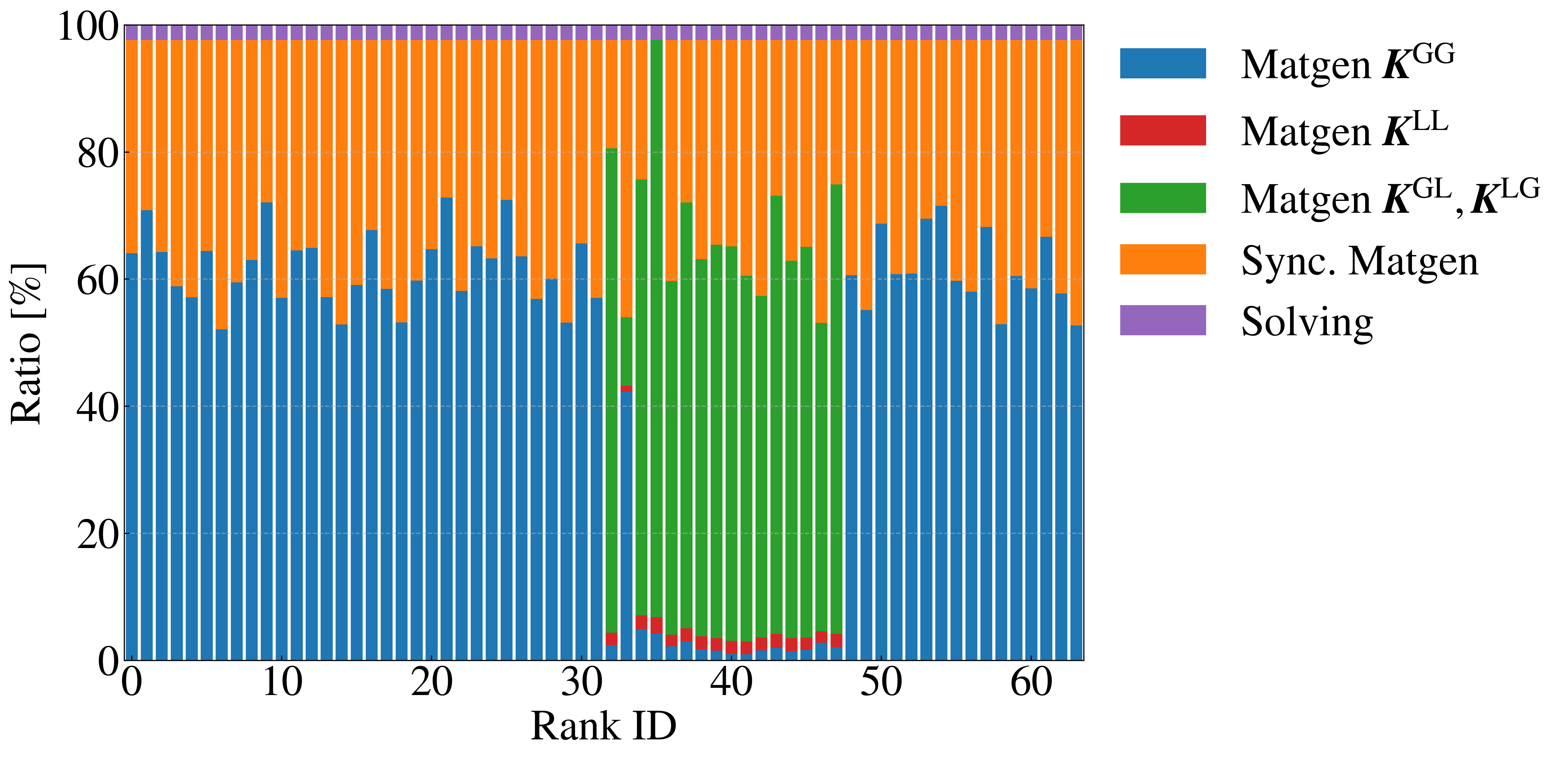}
        \par\small (b) $\alpha_3=20$
    \end{minipage}%
    \hfill
    \begin{minipage}{0.49\textwidth}
        \centering
        \includegraphics[bb=0 0 1094 530, width=\textwidth]{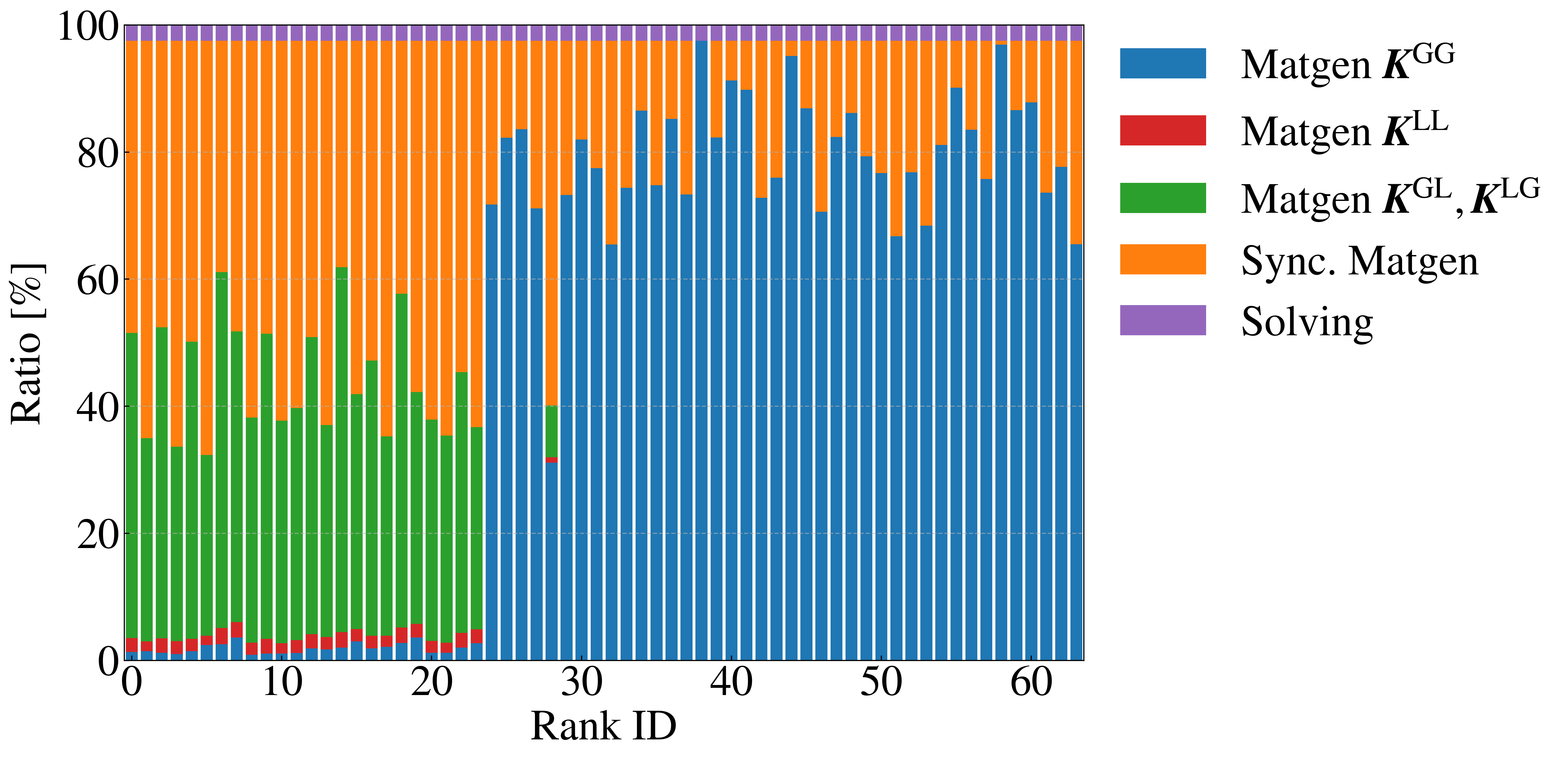}
        \par\small (c) $\alpha_3=40$
    \end{minipage}

    \vspace{0.3cm}

    \begin{minipage}{0.49\textwidth}
        \centering
        \includegraphics[bb=0 0 1094 530, width=\textwidth]{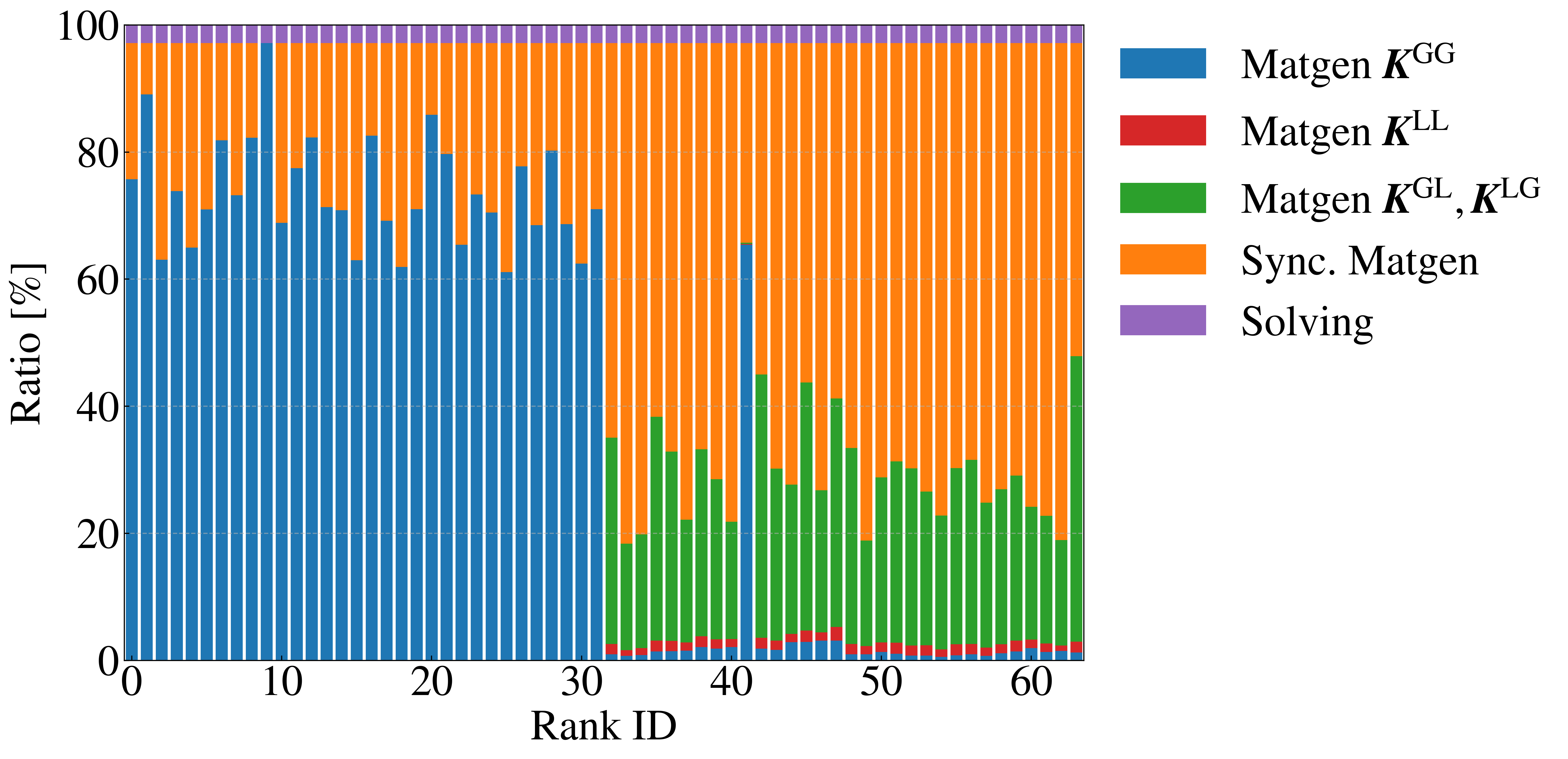}
        \par\small (d) $\alpha_3=60$
    \end{minipage}%
    \hfill
    \begin{minipage}{0.49\textwidth}
        \centering
        \includegraphics[bb=0 0 1094 530, width=\textwidth]{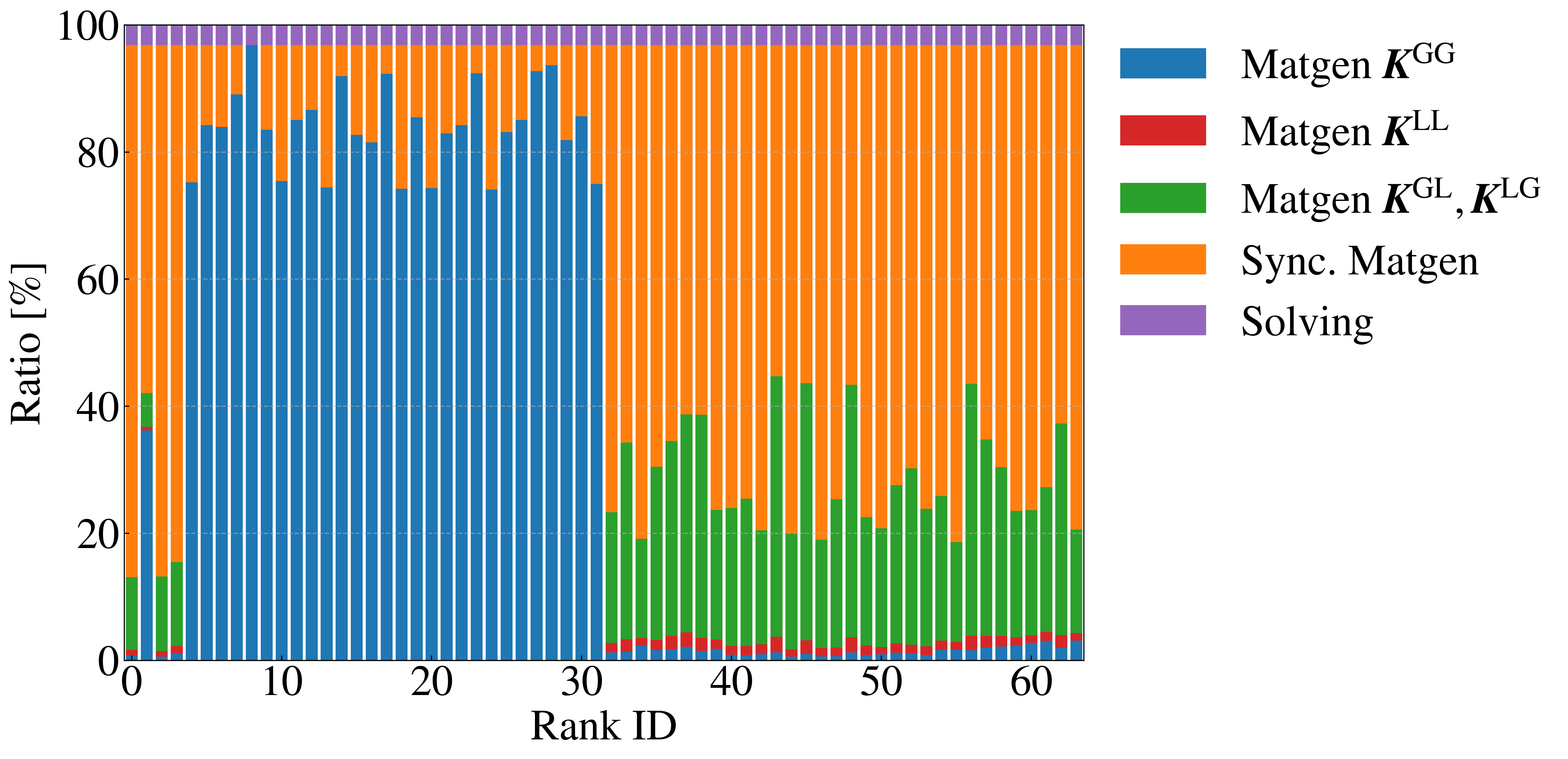}
        \par\small (e) $\alpha_3=80$
    \end{minipage}

    \caption{Breakdown of the measured time for each process under different weighting conditions at 64 processes. The blue, red, and green bars denote the generation times of $\boldsymbol{K}^{\mathrm{GG}}$, $\boldsymbol{K}^{\mathrm{LL}}$, and $\boldsymbol{K}^{\mathrm{GL}}/\boldsymbol{K}^{\mathrm{LG}}$, respectively. The orange bars denote the synchronization waiting time during the matrix generation phase, and the purple bars denote the linear solver time. Increasing $\alpha_3$ distributes the generation workload of $\boldsymbol{K}^{\mathrm{GL}}/\boldsymbol{K}^{\mathrm{LG}}$ over more subdomains, whereas the generation workload of $\boldsymbol{K}^{\mathrm{GG}}$ tends to become concentrated in fewer subdomains.}
    \label{fig:six_time_ratio}
\end{figure}

Next, to investigate how the introduction of node weights changes the workload distribution among subdomains, the breakdown of measured time for each process was compared at 64 processes.
\figurename~\ref{fig:six_time_ratio} shows the results for the no-weight case and for the different values of $\alpha_3$.
The horizontal axis represents the process number, and the vertical axis represents the fraction of measured time.

The figure shows that, as $\alpha_3$ increases, the generation workload of the coupling matrices $\boldsymbol{K}^{\mathrm{GL}}$ and $\boldsymbol{K}^{\mathrm{LG}}$ is distributed over a larger number of subdomains.
This occurs because increasing $\alpha_3$ assigns relatively larger weights to the nodes involved in coupling-matrix generation, making these nodes more likely to be distributed among multiple subdomains through graph partitioning.

On the other hand, when $\alpha_3$ becomes excessively large, the generation workload of $\boldsymbol{K}^{\mathrm{GL}}/\boldsymbol{K}^{\mathrm{LG}}$ is distributed more broadly, whereas the generation workload of $\boldsymbol{K}^{\mathrm{GG}}$ tends to become concentrated in a smaller number of subdomains.
This result indicates a trade-off between the workload distributions associated with the generation of $\boldsymbol{K}^{\mathrm{GL}}/\boldsymbol{K}^{\mathrm{LG}}$ and $\boldsymbol{K}^{\mathrm{GG}}$ in load balancing based on node weights.

These results confirm that incorporating matrix-generation costs into the node weights enables the localized computational workload associated with inter-mesh coupling to be distributed over multiple subdomains.

\begin{figure}
    \centering
    \includegraphics[bb=0 0 532 271, width=13cm]{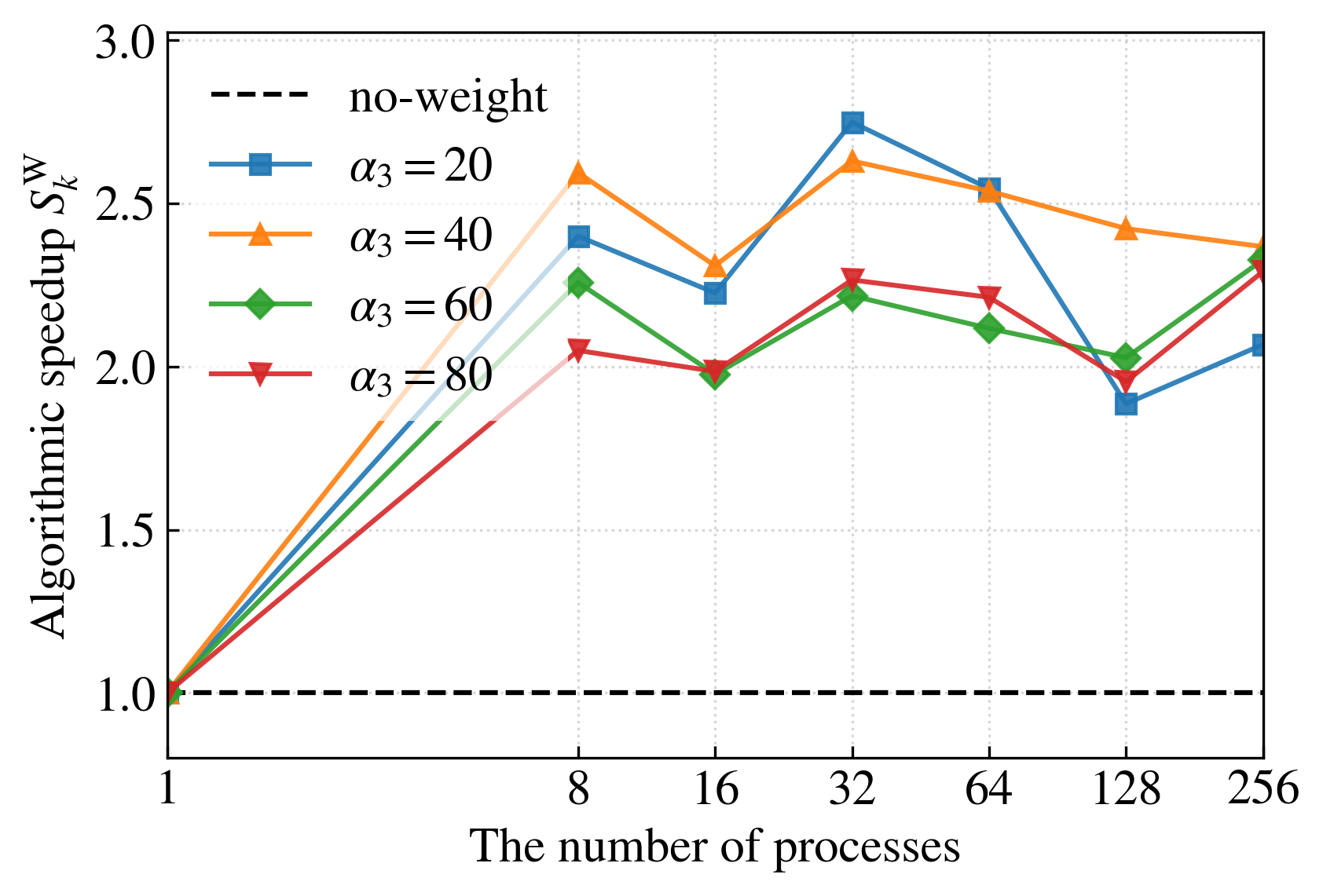}
    \caption{Weighting-induced speedup $S_k^{\mathrm{W}}$ relative to the no-weight baseline for different numbers of processes. The solid lines show the results for different weighting parameters: $\alpha_3=20$ (blue), $\alpha_3=40$ (orange), $\alpha_3=60$ (green), and $\alpha_3=80$ (red). The black dashed line denotes the no-weight baseline. A value greater than 1.0 indicates a reduction in total measured time due to the introduction of computational-cost-based node weights.}
    \label{fig:algorithmic_acceleration}
\end{figure}

Next, to evaluate the effect of the improved static load balancing described above on the total measured time, the weighting-induced speedup $S_k^{\mathrm{W}}$ relative to the no-weight case was evaluated.
\begin{align}
    S_k^{\mathrm{W}}
    =
    \frac{T_{k,\mathrm{noweight}}}{T_{k,\mathrm{weighted}}} ,
    \label{eq:weight_speedup}
\end{align}
where $T_{k,\mathrm{noweight}}$ denotes the total measured time for the no-weight case with $k$ processes, and $T_{k,\mathrm{weighted}}$ denotes the total measured time with $k$ processes when the weights $\alpha_1$, $\alpha_2$, and $\alpha_3$ are applied.
Thus, $S_k^{\mathrm{W}}>1$ indicates that the total measured time is reduced by weighting.

\figurename~\ref{fig:algorithmic_acceleration} shows the weighting-induced speedup $S_k^{\mathrm{W}}$.
The horizontal axis represents the number of parallel processes, and the vertical axis represents the speedup relative to the no-weight case.
The figure shows that, for all weighting conditions examined, the total measured time is reduced relative to the no-weight case.
This demonstrates that the improvement in static load balancing observed in \figurename~\ref{fig:stacked_bar_alpha3} and \figurename~\ref{fig:six_time_ratio} also leads to improved execution performance.

\begin{figure}
    \centering
    \includegraphics[bb=0 0 412 268, width=10cm]{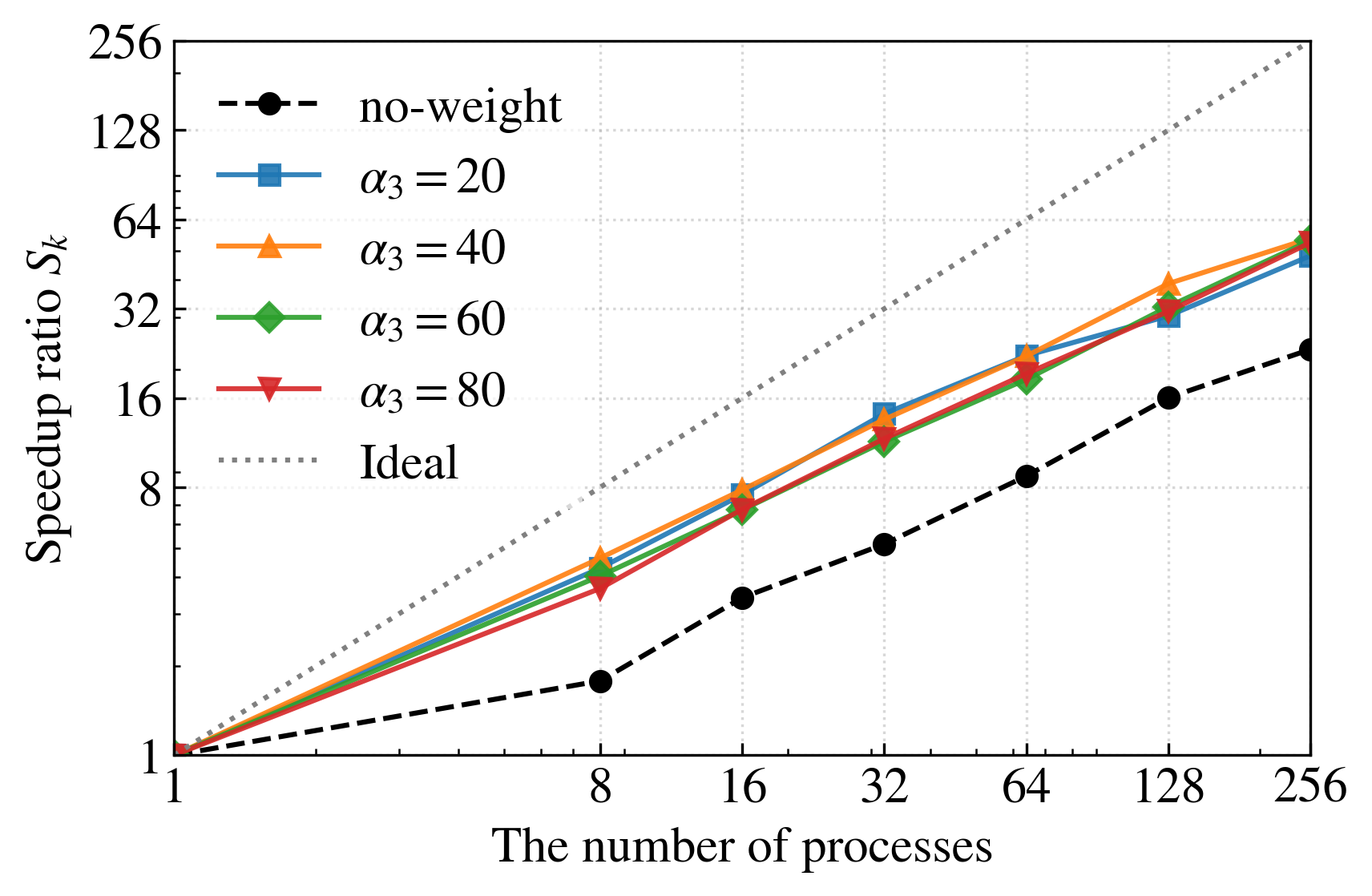}
    \caption{Strong-scaling performance in terms of the speedup $S_k$ based on the total measured time. The solid lines show the results for different weighting parameters: $\alpha_3=20$ (blue), $\alpha_3=40$ (orange), $\alpha_3=60$ (green), and $\alpha_3=80$ (red). The black dashed line denotes the no-weight baseline, and the gray dotted line denotes ideal linear scaling.}
    \label{fig:sfem_parastudy_speedup}
\end{figure}

\begin{figure}
    \centering
    \includegraphics[bb=0 0 412 268, width=10cm]{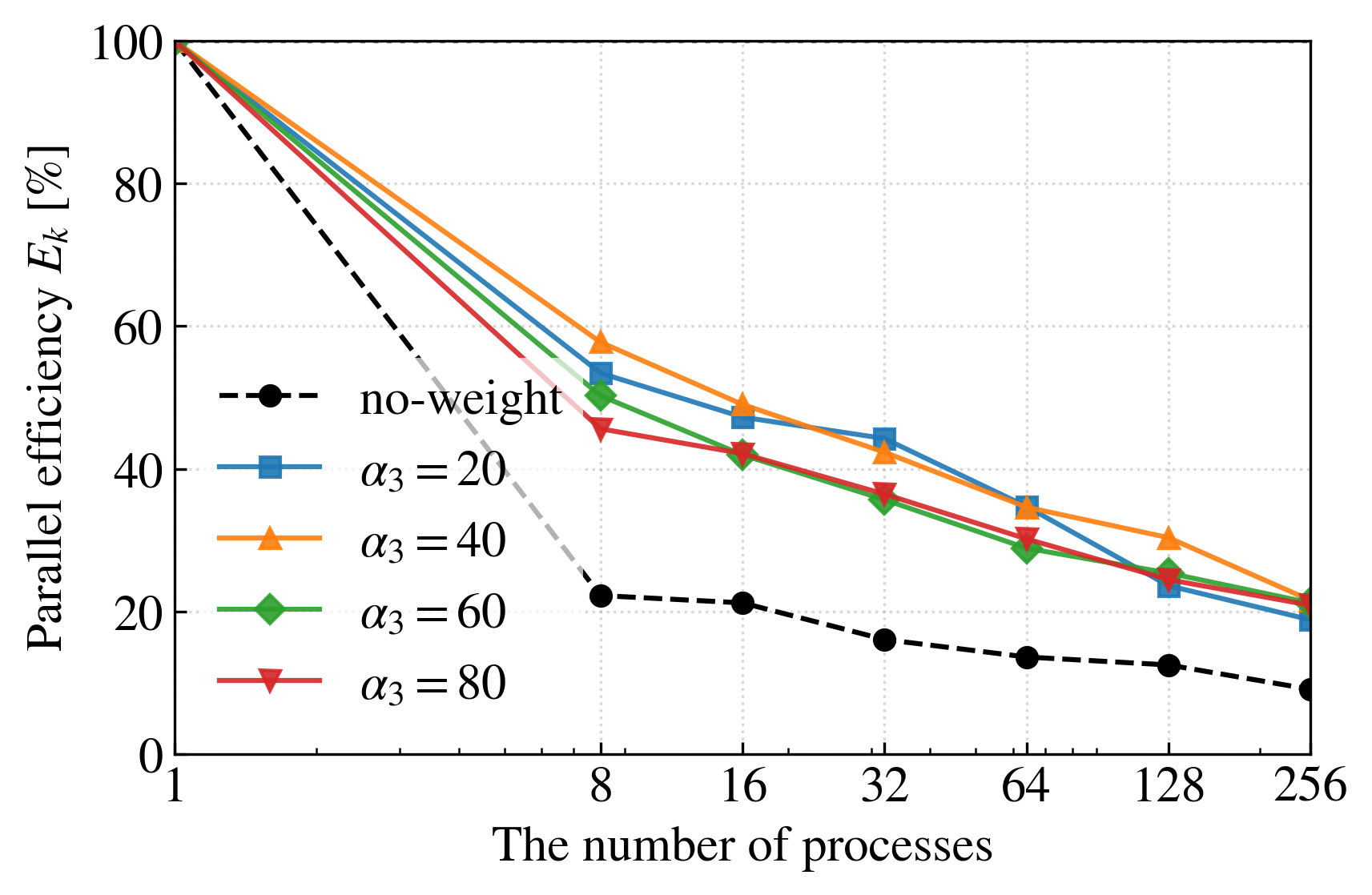}
    \caption{Strong-scaling performance in terms of the parallel efficiency $E_k$ based on the total measured time. The solid lines show the results for different weighting parameters: $\alpha_3=20$ (blue), $\alpha_3=40$ (orange), $\alpha_3=60$ (green), and $\alpha_3=80$ (red). The black dashed line denotes the no-weight baseline, and the gray dotted line denotes the ideal parallel efficiency.}
    \label{fig:sfem_parastudy_efficiency}
\end{figure}

Finally, the effect of introducing node weights on strong-scaling performance was evaluated.
The speedup $S_k$ and parallel efficiency $E_k$ defined in Eqs.~\eqref{eq:s_k} and \eqref{eq:e_k}, respectively, were used.
The evaluation was based on the total measured time.

The speedup $S_k$ is shown in \figurename~\ref{fig:sfem_parastudy_speedup}, and the parallel efficiency $E_k$ is shown in \figurename~\ref{fig:sfem_parastudy_efficiency}.
In both figures, the horizontal axis represents the number of parallel processes, whereas the vertical axes represent the speedup and parallel efficiency, respectively.
The black dashed lines show the results for the no-weight case, and the gray dotted lines indicate ideal strong scaling.

These results show that introducing node weights that account for the coupling-matrix generation workload substantially improves both the speedup and parallel efficiency compared with the no-weight case.
For example, at 256 processes, the speedup was approximately 23 in the no-weight case, whereas it increased to approximately 55 for $\alpha_3=40$.

Overall, the results of this subsection demonstrate that incorporating the matrix-generation costs of BSFEM into the node weights of the proposed graph enables the nonuniform computational workload arising from inter-mesh coupling to be reflected in the domain decomposition, thereby improving static load balancing.
Consequently, the synchronization waiting time is reduced, and both the total measured time and parallel performance are improved.
This verification is not intended to provide a method for determining optimal node weights; rather, it demonstrates the utility of the proposed method by showing that additional information on computational workload can be incorporated into the same generalized graph.

\section{Conclusions}
\label{sec:conclusions}
In this study, we proposed a domain-decomposed parallelization framework for BSFEM. The proposed method represents interactions among computational points within and across the superimposed B-spline and Lagrange meshes in a single generalized graph.
The resulting graph partition assigns degrees of freedom and elements to processes. It also defines the MPI communication structure.
To the authors' knowledge, this framework realizes the first domain-decomposition-based distributed-memory parallel computation for an SFEM-based method.

Verification using a three-dimensional steady Poisson problem confirmed that the accuracy of the BSFEM solution was maintained as the number of subdomains was varied. Strong-scaling tests showed that the spatially nonuniform workload associated with coupling-matrix generation causes load imbalance and synchronization waiting time, thereby limiting parallel performance. As a demonstration of the utility of the generalized graph representation, matrix-generation costs were incorporated into the node weights. The resulting cost-weighted graph partitioning improved static load balancing, reduced synchronization waiting time and total measured time, and improved parallel performance.

The present study was limited to BSFEM with fixed meshes. Future work will include automatic determination of node weights, extension to nonlinear and unsteady problems, and generalization of the proposed framework to other multiple-mesh approaches, including conventional SFEM formulations and overset methods.

\section*{Acknowledgments}

This work was supported by JST FOREST under Grant No.~JPMJFR215S,
JST ACT-X under Grant No.~JPMJAX24LG,
the Joint Usage/Research Center for Interdisciplinary Large-scale
Information Infrastructures (JHPCN) under Grant No.~jh240017,
and JSPS KAKENHI under Grant Nos.~23H00475, 24KJ0499, and 26K02921.

\section*{Data availability}

The data that support the findings of this study are available from
the corresponding author upon reasonable request.

\printcredits

\section*{Declaration of competing interest}

The authors declare that they have no known competing financial
interests or personal relationships that could have appeared to
influence the work reported in this paper.


\appendix

\counterwithin{equation}{section}
\counterwithin{figure}{section}
\counterwithin{table}{section}

\section{Computational complexity of submatrix generation}
\label{sec:appendix_A}

This section summarizes the computational complexity required to
generate each submatrix in BSFEM.
In the following,
$n_{\mathrm{el}}^{\mathrm{G}}$ and
$n_{\mathrm{el}}^{\mathrm{L}}$ denote the numbers of elements in the
global and local meshes, respectively;
$n_{\mathrm{ip}}^{\mathrm{G}}$ and
$n_{\mathrm{ip}}^{\mathrm{L}}$ denote the numbers of integration
points per element;
$n_{\mathrm{ldof}}^{\mathrm{G}}$ and
$n_{\mathrm{ldof}}^{\mathrm{L}}$ denote the numbers of nodes in an
element in the global and local meshes, respectively; and $d$ denotes the number of degrees of freedom per node.
For the scalar Poisson equation considered in this study, \(d=1\).

The submatrix $\boldsymbol{K}^{\mathrm{GG}}$ associated with the global
mesh is generated using standard finite element procedures with
B-spline basis functions.
Its computational complexity is given by
\begin{equation}
    \mathcal{O}
    \left(
    n_{\mathrm{el}}^{\mathrm{G}}
    n_{\mathrm{ip}}^{\mathrm{G}}
    d^2
    \left(n_{\mathrm{ldof}}^{\mathrm{G}}\right)^2
    \right).
    \label{eq:appendix_KGG_complexity}
\end{equation}

The submatrix $\boldsymbol{K}^{\mathrm{LL}}$ associated with the local
mesh is generated using standard finite element procedures with
Lagrange basis functions.
Its computational complexity is given by
\begin{equation}
    \mathcal{O}
    \left(
    n_{\mathrm{el}}^{\mathrm{L}}
    n_{\mathrm{ip}}^{\mathrm{L}}
    d^2
    \left(n_{\mathrm{ldof}}^{\mathrm{L}}\right)^2
    \right).
    \label{eq:appendix_KLL_complexity}
\end{equation}

The coupling matrices $\boldsymbol{K}^{\mathrm{GL}}$ and
$\boldsymbol{K}^{\mathrm{LG}}$ are generated by evaluating the
interactions between the local basis functions and the corresponding
global basis functions at the integration points of the local elements.
Because the corresponding global elements and basis function values
are stored in the coupling map, no geometric search is required during
matrix generation.
Therefore, the computational complexity of generating the coupling
matrices is
\begin{equation}
    \mathcal{O}
    \left(
    n_{\mathrm{el}}^{\mathrm{L}}
    n_{\mathrm{ip}}^{\mathrm{L}}
    d^2
    n_{\mathrm{ldof}}^{\mathrm{G}}
    n_{\mathrm{ldof}}^{\mathrm{L}}
    \right).
    \label{eq:appendix_coupling_complexity}
\end{equation}

\section{Formulation of node weights in BSFEM}
\label{sec:appendix_B}

In the previous section, the computational complexity required to
generate each submatrix in BSFEM was summarized.
In this section, this complexity estimate is reduced to an element-wise
local matrix-generation cost.
The cost is then distributed among the associated nodes to explicitly
define the node weight $w_i^{\mathrm{V}}$ introduced in
Section~\ref{sq:weighted_graph}.

In the global matrix generation of BSFEM, the types of elements involved
and the combinations of basis functions differ among
$\boldsymbol{K}^{\mathrm{GG}}$, $\boldsymbol{K}^{\mathrm{LL}}$, and the
coupling matrices $\boldsymbol{K}^{\mathrm{GL}}$ and
$\boldsymbol{K}^{\mathrm{LG}}$.
Therefore, in the formulation of the node weights, the local
computational costs associated with these submatrices are treated
separately.

In this study, the element-wise matrix-generation cost $w_e$ is defined
based on the computational complexity of the submatrices described in
the previous section.
The node weight $w_i^{\mathrm{V}}$ is then obtained by distributing this
cost among the associated nodes and accumulating the contributions at
each node.

Here, the geometric database and coupling map are assumed to have
already been constructed.
The estimate therefore focuses on the dominant operation count
associated with the evaluation of element matrices for each element.
Consequently, the geometric search cost required for constructing the
coupling map is not included in the following complexity estimates.
In the following, $n_{\mathrm{ip}}$ denotes the number of integration
points in an element, $n_{\mathrm{ldof}}$ denotes the number of nodes
in an element, and $d$ denotes the number of degrees of freedom per
node.
For the scalar Poisson problem considered in this study, $d=1$.

\begin{itemize}

\item The local computational cost for
$\boldsymbol{K}^{\mathrm{GG}}$ on a global element $e$ is defined as
\begin{equation}
    w_e^{\mathrm{GG}}
    =
    \alpha_1
    n_{\mathrm{ip}}^{\mathrm{G}}
    d^2
    \left(n_{\mathrm{ldof}}^{\mathrm{G}}\right)^2.
    \label{eq:appendix_weight_KGG}
\end{equation}

\item The local computational cost for
$\boldsymbol{K}^{\mathrm{LL}}$ on a local element $e$ is defined as
\begin{equation}
    w_e^{\mathrm{LL}}
    =
    \alpha_2
    n_{\mathrm{ip}}^{\mathrm{L}}
    d^2
    \left(n_{\mathrm{ldof}}^{\mathrm{L}}\right)^2.
    \label{eq:appendix_weight_KLL}
\end{equation}

\item The coupling matrices $\boldsymbol{K}^{\mathrm{GL}}$ and
$\boldsymbol{K}^{\mathrm{LG}}$ are generated between a local element
and the global elements that spatially overlap it.
Here, $w_e^{\mathrm{GL,LG}}$ denotes the total cost of the coupling
terms associated with the generation of
$\boldsymbol{K}^{\mathrm{GL}}$ and
$\boldsymbol{K}^{\mathrm{LG}}$.
\begin{equation}
    w_e^{\mathrm{GL,LG}}
    =
    \alpha_3
    n_{\mathrm{ip}}^{\mathrm{L}}
    d^2
    n_{\mathrm{ldof}}^{\mathrm{L}}
    n_{\mathrm{ldof}}^{\mathrm{G}}.
    \label{eq:appendix_weight_coupling}
\end{equation}

\end{itemize}

Here, $\alpha_1$, $\alpha_2$, and $\alpha_3$ are coefficients that
convert the theoretical operation counts into the actual computational
load.

To convert these element-based costs into the node weight
$w_i^{\mathrm{V}}$, the costs $w_e^{\mathrm{GG}}$ and
$w_e^{\mathrm{LL}}$ are uniformly distributed among the constituent
nodes $V_e$ of the corresponding element $e$.
In contrast, for the coupling matrix cost $w_e^{\mathrm{GL,LG}}$, the
cost is equally divided between the Lagrange nodes and the B-spline
nodes.
This distribution is not an exact decomposition of the operation count,
but rather an approximate weighting strategy for load balancing.

The final node weight $w_i^{\mathrm{V}}$ is defined by the following
summations, depending on the type of node.
Let $\mathcal{E}_i^{\mathrm{G}}$ and
$\mathcal{E}_i^{\mathrm{L}}$ denote the sets of global and local
elements that contain node $v_i$, respectively.
Let $\mathcal{E}^{\mathrm{L}}$ denote the set of all elements in a local
mesh.
Furthermore, for a local element $e \in \mathcal{E}^{\mathrm{L}}$, let
$\mathcal{M}_e^{\mathrm{G}}$ denote the set of global elements crossed
by the integration points in $e$, and define its cardinality as
$m_e^{\mathrm{G}} = |\mathcal{M}_e^{\mathrm{G}}|$.

First, for a B-spline node $v_i \in V^{\mathrm{G}}$, the weight is
given by the sum of the contribution from the global elements to which
the node belongs and the contribution from coupling-matrix generation
associated with the overlapping local elements:
\begin{equation}
w_i^{\mathrm{V}}
= 
\sum_{e \in \mathcal{E}_i^{\mathrm{G}}}
\frac{w_e^{\mathrm{GG}}}
{n_{\mathrm{ldof}}^{\mathrm{G}}}
+
\sum_{e \in \mathcal{E}^{\mathrm{L}}}
\sum_{e^{\mathrm{G}} \in \mathcal{M}_e^{\mathrm{G}}}
\frac{0.5\,w_e^{\mathrm{GL,LG}}}
{m_e^{\mathrm{G}}
n_{\mathrm{ldof}}^{\mathrm{G}}}
\delta_{i,e^{\mathrm{G}}}.
\label{eq:appendix_global_node_weight}
\end{equation}

Here, $\delta_{i,e^{\mathrm{G}}}$ is an indicator function that takes
the value $1$ if node $v_i$ belongs to the global element
$e^{\mathrm{G}}$, and $0$ otherwise.
The second term on the right-hand side represents the operation of
uniformly distributing half of the coupling-matrix generation cost
generated by the local element $e$ among all nodes constituting the
associated $m_e^{\mathrm{G}}$ global elements.
For simplicity, the coupling-matrix generation cost arising from the
local element $e$ is uniformly distributed among the corresponding
$m_e^{\mathrm{G}}$ global elements.

Next, for a Lagrange node $v_i \in V^{\mathrm{L}}$, the weight is given
by the sum of the contribution from $\boldsymbol{K}^{\mathrm{LL}}$
associated with the local elements to which the node belongs and the
contribution from coupling-matrix generation arising from those local
elements.
Under the assumption
$\Omega^{\mathrm{L}} \subseteq \Omega^{\mathrm{G}}$, the coupling terms
are evaluated for all local elements containing node $v_i$.
Thus, the weight is expressed as
\begin{equation}
w_i^{\mathrm{V}}
=
\sum_{e \in \mathcal{E}_i^{\mathrm{L}}}
\frac{
w_e^{\mathrm{LL}}
+ 0.5\,w_e^{\mathrm{GL,LG}}
}{
n_{\mathrm{ldof}}^{\mathrm{L}}
}.
\label{eq:appendix_local_node_weight}
\end{equation}

It should be noted that the optimal ratio among the coefficients
$\alpha_1$, $\alpha_2$, and $\alpha_3$ is difficult to determine
a priori.
This is because the computational cost of each term is affected not
only by the number of arithmetic operations, but also by
implementation-level factors such as memory access locality and the
frequency of indirect references.
Therefore, in this study, these coefficients are treated as tuning
parameters depending on the computational environment, and their
optimal values are identified through numerical experiments on the
target machine.

\bibliographystyle{cas-model2-names}

\bibliography{cas-refs}

@article{thames1977numerical,
  title={Numerical solutions for viscous and potential flow about arbitrary two-dimensional bodies using body-fitted coordinate systems},
  author={Thames, Frank C and Thompson, Joe F and Mastin, C Wayne and Walker, Ray L},
  journal={Journal of Computational Physics},
  volume={24},
  number={3},
  pages={245--273},
  year={1977},
  publisher={Elsevier}
}

@article{hughes2005isogeometric,
  title={Isogeometric analysis: {CAD}, finite elements, {NURBS}, exact geometry and mesh refinement},
  author={Hughes, Thomas JR and Cottrell, John A and Bazilevs, Yuri},
  journal={Computer methods in applied mechanics and engineering},
  volume={194},
  number={39--41},
  pages={4135--4195},
  year={2005},
  publisher={Elsevier}
}

@article{el2022anisotropic,
  title={Anisotropic adaptive body-fitted meshes for {CFD}},
  author={El Aouad, Sacha and Larcher, Aur{\'e}lien and Hachem, Elie},
  journal={Computer Methods in Applied Mechanics and Engineering},
  volume={400},
  pages={115562},
  year={2022},
  publisher={Elsevier}
}

@article{mittal2005immersed,
  title={Immersed boundary methods},
  author={Mittal, Rajat and Iaccarino, Gianluca},
  journal={Annual Review of Fluid Mechanics},
  volume={37},
  number={1},
  pages={239--261},
  year={2005},
  publisher={Annual Reviews}
}

@inproceedings{benek1983flexible,
  title={A flexible grid embedding technique with application to the {Euler} equations},
  author={Benek, J and Steger, J and Dougherty, F Carroll},
  booktitle={6th computational fluid dynamics conference Danvers},
  pages={1944},
  year={1983}
}

@article{steger1987use,
  title={On the use of composite grid schemes in computational aerodynamics},
  author={Steger, Joseph L and Benek, John A},
  journal={Computer Methods in Applied Mechanics and Engineering},
  volume={64},
  number={1--3},
  pages={301--320},
  year={1987},
  publisher={Elsevier}
}

@article{chesshire1990composite,
  title={Composite overlapping meshes for the solution of partial differential equations},
  author={Chesshire, G and Henshaw, William D},
  journal={Journal of Computational Physics},
  volume={90},
  number={1},
  pages={1--64},
  year={1990},
  publisher={Elsevier}
}

@article{tang2003overset,
  title={An overset-grid method for {3D} unsteady incompressible flows},
  author={Tang, HS and Jones, S Casey and Sotiropoulos, Fotis},
  journal={Journal of Computational Physics},
  volume={191},
  number={2},
  pages={567--600},
  year={2003},
  publisher={Elsevier}
}

@article{appelo2012numerical,
  title={Numerical methods for solid mechanics on overlapping grids: Linear elasticity},
  author={Appel{\"o}, Daniel and Banks, Jeffrey W and Henshaw, William D and Schwendeman, Donald W},
  journal={Journal of Computational Physics},
  volume={231},
  number={18},
  pages={6012--6050},
  year={2012},
  publisher={Elsevier}
}

@article{henshaw2009composite,
  title={A composite grid solver for conjugate heat transfer in fluid--structure systems},
  author={Henshaw, William D and Chand, Kyle K},
  journal={Journal of Computational Physics},
  volume={228},
  number={10},
  pages={3708--3741},
  year={2009},
  publisher={Elsevier}
}

@article{fish1992s,
  title={The s-version of the finite element method},
  author={Fish, Jacob},
  journal={Computers \& Structures},
  volume={43},
  number={3},
  pages={539--547},
  year={1992},
  publisher={Elsevier}
}

@article{fish1993multiscale,
  title={Multiscale finite element method for a locally nonperiodic heterogeneous medium},
  author={Fish, Jacob and Wagiman, Amir},
  journal={Computational mechanics},
  volume={12},
  number={3},
  pages={164--180},
  year={1993},
  publisher={Springer}
}

@article{FISH1994135,
title = {On adaptive multilevel superposition of finite element meshes for linear elastostatics},
journal = {Applied Numerical Mathematics},
volume = {14},
number = {1},
pages = {135--164},
year = {1994},
issn = {0168-9274},
author = {J. Fish and S. Markolefas and R. Guttal and P. Nayak}
}

@article{lee2004combined,
  title={Combined extended and superimposed finite element method for cracks},
  author={Lee, Sang-Ho and Song, Jeong-Hoon and Yoon, Young-Cheol and Zi, Goangseup and Belytschko, Ted},
  journal={International Journal for Numerical Methods in Engineering},
  volume={59},
  number={8},
  pages={1119--1136},
  year={2004},
  publisher={Wiley Online Library}
}

@article{okada2005fracture,
  title={On fracture analysis using an element overlay technique},
  author={Okada, Hiroshi and Endoh, Sayaka and Kikuchi, Masanori},
  journal={Engineering fracture mechanics},
  volume={72},
  number={5},
  pages={773--789},
  year={2005},
  publisher={Elsevier}
}

@article{nakasumi2008crack,
  title={Crack growth analysis using mesh superposition technique and {X-FEM}},
  author={Nakasumi, Shogo and Suzuki, Katsuyuki and Ohtsubo, Hideomi},
  journal={International journal for numerical methods in engineering},
  volume={75},
  number={3},
  pages={291--304},
  year={2008},
  publisher={Wiley Online Library}
}

@article{kikuchi2014fatigue,
  title={Fatigue crack growth simulation in heterogeneous material using s-version {FEM}},
  author={Kikuchi, Masanori and Wada, Yoshitaka and Shintaku, Yuichi and Suga, Kazuhiro and Li, Yulong},
  journal={International Journal of Fatigue},
  volume={58},
  pages={47--55},
  year={2014},
  publisher={Elsevier}
}

@article{xu2018study,
  title={A study on the {S}-version {FEM} for a dynamic damage model},
  author={Xu, Qiang and Chen, Jianyun and Yue, Hongyuan and Li, Jing},
  journal={International Journal for Numerical Methods in Engineering},
  volume={115},
  number={4},
  pages={427--444},
  year={2018},
  publisher={Wiley Online Library}
}

@article{kishi2020dynamic,
  title={Dynamic crack propagation analysis based on the s-version of the finite element method},
  author={Kishi, Kota and Takeoka, Yuuki and Fukui, Tsutomu and Matsumoto, Toshiyuki and Suzuki, Katsuyuki and Shibanuma, Kazuki},
  journal={Computer Methods in Applied Mechanics and Engineering},
  volume={366},
  pages={113091},
  year={2020},
  publisher={Elsevier}
}

@article{sun2018variant,
  title={A variant of the s-version of the finite element method for concurrent multiscale coupling},
  author={Sun, Wei and Fish, Jacob and Dhia, Hachmi Ben},
  journal={International Journal for Multiscale Computational Engineering},
  volume={16},
  number={2},
  year={2018},
  publisher={Begel House Inc.}
}

@article{sakata2022mesh,
  title={Mesh superposition-based multiscale stress analysis of composites using homogenization theory and re-localization technique considering fiber location variation},
  author={Sakata, Sei-ichiro and Tanimasu, Shin},
  journal={International Journal for Numerical Methods in Engineering},
  volume={123},
  number={2},
  pages={505--529},
  year={2022},
  publisher={Wiley Online Library}
}

@article{wang2006moving,
  title={A moving superimposed finite element method for structural topology optimization},
  author={Wang, Shengyin and Wang, Michael Y},
  journal={International Journal for Numerical Methods in Engineering},
  volume={65},
  number={11},
  pages={1892--1922},
  year={2006},
  publisher={Wiley Online Library}
}

@article{yue2005adaptive,
  title={Adaptive superposition of finite element meshes in elastodynamic problems},
  author={Yue, Z and Robbins Jr, DH},
  journal={International journal for numerical methods in engineering},
  volume={63},
  number={11},
  pages={1604--1635},
  year={2005},
  publisher={Wiley Online Library}
}

@article{FISH1993363,
title = {Adaptive s-method for linear elastostatics},
journal = {Computer Methods in Applied Mechanics and Engineering},
volume = {104},
number = {3},
pages = {363--396},
year = {1993},
issn = {0045-7825},
author = {Jacob Fish and Stilianos Markolefas}
}

@article{he2023strategy,
  title={Strategy for accurately and efficiently modelling an internal traction-free boundary based on the s-version finite element method: Problem clarification and solutions verification},
  author={He, Tianyu and Mitsume, Naoto and Yasui, Fumitaka and Morita, Naoki and Fukui, Tsutomu and Shibanuma, Kazuki},
  journal={Computer Methods in Applied Mechanics and Engineering},
  volume={404},
  pages={115843},
  year={2023},
  publisher={Elsevier}
}

@article{fan2008rs,
  title={The rs-method for material failure simulations},
  author={Fan, R and Fish, Jacob},
  journal={International journal for numerical methods in engineering},
  volume={73},
  number={11},
  pages={1607--1623},
  year={2008},
  publisher={Wiley Online Library}
}

@article{ooya2009linear,
  title={On the linear dependencies of interpolation functions in s-version finite element method},
  author={Ooya, Takanori and Tanaka, Satoyuki and Okada, Hiroshi},
  journal={Journal of Computational Science and Technology},
  volume={3},
  number={1},
  pages={124--135},
  year={2009},
  publisher={The Japan Society of Mechanical Engineers}
}

@article{magome2024higher,
  title={Higher-continuity s-version of finite element method with {B}-spline functions},
  author={Magome, Nozomi and Morita, Naoki and Kaneko, Shigeki and Mitsume, Naoto},
  journal={Journal of Computational Physics},
  volume={497},
  pages={112593},
  year={2024},
  publisher={Elsevier}
}

@article{suwa2023parallel,
  title={Parallel parametric analysis approach based on an s-version {FEM} for fracture mechanics analysis in designing hole positions},
  author={Suwa, Hiroki and Yusa, Yasunori and Kuboki, Takashi},
  journal={Mechanical Engineering Journal},
  volume={10},
  number={3},
  pages={22--00462},
  year={2023},
  publisher={The Japan Society of Mechanical Engineers}
}

@article{karypis1998fast,
  title={A fast and high quality multilevel scheme for partitioning irregular graphs},
  author={Karypis, George and Kumar, Vipin},
  journal={SIAM Journal on scientific Computing},
  volume={20},
  number={1},
  pages={359--392},
  year={1998},
  publisher={SIAM}
}

@article{otoguro2017space,
  title={Space--time {VMS} computational flow analysis with isogeometric discretization and a general-purpose {NURBS} mesh generation method},
  author={Otoguro, Yuto and Takizawa, Kenji and Tezduyar, Tayfun E},
  journal={Computers \& Fluids},
  volume={158},
  pages={189--200},
  year={2017},
}

@article{Borden2011,
  title={Isogeometric finite element data structures based on {B{\'e}zier} extraction of {NURBS}},
  author={Borden, Michael J and Scott, Michael A and Evans, John A and Hughes, Thomas JR},
  journal={International Journal for Numerical Methods in Engineering},
  volume={87},
  number={1--5},
  pages={15--47},
  year={2011},
}

@misc{squid,
  author       = {{Cybermedia Center, Osaka University}},
  title        = {{Supercomputer for Quest to Unsolved Interdisciplinary Datascience (SQUID)}},
  howpublished = {\url{http://www.hpc.cmc.osaka-u.ac.jp/squid/}},
  note         = {Accessed: February 26, 2026}
}

@misc{monolis,
  author       = {Morita, Naoki and Mitsume, Naoto},
  title        = {{Metagraph-Oriented Network for Linear Iterative Solvers (MONOLIS)}},
  howpublished = {\url{https://www.kz.tsukuba.ac.jp/~nmorita/monolis.html}},
  year         = {2026},
  note         = {Accessed: February 26, 2026}
}

@book{roache1998verification,
  title     = {Verification and Validation in Computational Science and Engineering},
  author    = {Roache, Patrick J.},
  volume    = {895},
  year      = {1998},
  publisher = {Hermosa Publishers},
  address   = {Albuquerque}
}



\end{document}